\documentclass[]{article}

\usepackage{mathrsfs}
\usepackage{amssymb}
\usepackage{amsfonts}
\usepackage{graphicx}
\usepackage{subfigure}
\usepackage{float}
\usepackage[titletoc]{appendix}
\usepackage{array}
\usepackage{amsmath}
\usepackage{epsf}
\usepackage{indentfirst}
\usepackage{amsthm}
\usepackage{color}
\usepackage{hyperref}
\usepackage{titletoc}

\numberwithin{equation}{section}
\newtheorem{theorem}{Theorem}[section]

\newtheorem{remark}[theorem]{Remark}

\newcommand{\rd}{\mathrm{d}}

\begin{document}
\title{Asymptotic-preserving exponential methods for the quantum 
Boltzmann equation with high-order accuracy\footnote{This work was partially 
supported by RNMS11-07444 (KI-Net) and by PRIN 2009 project ``Advanced numerical methods for kinetic equations and balance laws with source terms''.}}
\author{Jingwei Hu\footnote{Institute for Computational Engineering and 
Sciences (ICES), The University of Texas at Austin, 201 East 24th St, Stop 
C0200, Austin, TX 78712, USA (hu@ices.utexas.edu).},\ \  Qin Li\footnote{Department 
of Computing + Mathematical Sciences (CMS), The Annenberg Center, California Institute of 
Technology, 1200 E. California Blvd., Pasadena, CA 91125, USA (qinli@caltech.edu).}, \ and Lorenzo 
Pareschi\footnote{Department of Mathematics and Computer Science, University of Ferrara, Via Machiavelli 35, 44121 Ferrara, Italy (lorenzo.pareschi@unife.it)}}
\maketitle

\begin{abstract}
In this paper we develop high order asymptotic preserving methods for the 
spatially inhomogeneous quantum Boltzmann equation. We follow the work in Li 
and Pareschi \cite{LP_ExpRKinhomoBoltzmann} where asymptotic preserving 
exponential Runge-Kutta methods for the classical inhomogeneous Boltzmann 
equation were constructed. A major difficulty here is related to the non 
Gaussian steady states characterizing the quantum kinetic behavior. We show 
that the proposed schemes work with high-order accuracy uniformly in time for 
all Planck constants ranging from classical regime to quantum regime, and all 
Knudsen number ranging from kinetic regime to fluid regime. Computational 
results are presented for both Bose gas and Fermi gas.
\end{abstract}

{\small 
{\bf Key words.} Quantum Boltzmann equation, asymptotic preserving methods, exponential Runge-Kutta schemes.

{\bf AMS subject classifications.} 65L04, 65L06, 35Q20, 82C10.
}

\section{Introduction}

The quantum Boltzmann equation (QBE), also known as the Nordheim-Uehling-Uhlenbeck equation, describes the nonequilibrium dynamics of a dilute quantum gas consisting of elementary particles of bosons or fermions  \cite{CC}. By including quantum mechanical effects in the collisional process, the equation models a wider range of particle behaviors than the usual Boltzmann equation of classical particles. This is because the latter can be treated as a sub-model under a certain classical limit (Planck constant approaching zero). The QBE and its variants have many applications in science and engineering, including the kinetic description of Bose-Einstein condensate \cite{ST95,Spohn10}, and the modeling of electron interactions in semiconductor devices \cite{MRS, Jungel}.
 
In this paper we design a class of high order numerical methods for quantum Boltzmann 
equation, that is accurate and efficient in both kinetic and hydrodynamic regimes for all 
Planck constants. In kinetic theory, the time discretization represents a 
computational challenge in the construction of numerical methods, especially 
in stiff regimes, when the collisional scale becomes dominant over the 
transport of particles, and the fluid-dynamic limit is achieved. To resolve 
the collision term, the time step is severely controlled by the Knudsen number 
for numerical stability if explicit schemes are used. On the other hand, the 
use of implicit schemes allows larger time steps but presents considerable limitations 
in most applications since the collision operator is usually highly nonlinear 
and nonlocal.

Many techniques have been developed to address such issues in recent years, 
and we specifically mention the micro-macro 
decomposition~\cite{BLM08}, the BGK penalization 
method~\cite{FJ_APGeneral}, and the exponential Runge-Kutta
methods~\cite{DP_ExpRK,GPT_RelaxationNonlinearKinetic}. The feature shared among these techniques is that 
the schemes are unconditionally stable, capturing the asymptotic limits 
automatically without time being resolved, and are numerically less complicated than other possible 
approaches, for example, the domain decomposition strategies and hybrid 
methods at different levels~\cite{PareschiCaflisch_implicitMC,TiwariKlar_adaptiveDD,DJM_KineticHydrodynamic,DimarcoPareschi_hybrid}. For a nice survey on asymptotic-preserving (AP) scheme for various kinds of systems see, for instance, the review paper by Jin~\cite{Jin_Review}. In the case of Boltzmann-type kinetic equations we refer to a recent review by Pareschi and Russo~\cite{PareschiRusso_reviewBoltzmannAP}.

In this work we extend the asymptotic preserving exponential Runge-Kutta method developed 
in~\cite{DP_ExpRK,LP_ExpRKinhomoBoltzmann} to the quantum Boltzmann equation. 
The extension to the multi-species Boltzmann equation could be found 
in~\cite{LY_ExpRKmultiBoltzmann}. We refer the reader 
to~\cite{HochOstermann_ExponentialIntegrator} for an introduction to time 
integration exponential techniques. New difficulties in the quantum case 
would be: 
\begin{itemize}
\item The steady states are not classical Maxwellian (Gaussian distribution) and to obtain the local equilibrium --- the quantum Maxwellian (Bose-Einstein or Fermi-Dirac distribution), a nonlinear system needs to be inverted; 
\item The methods developed need 
to be uniformly high order and efficient for all Planck constants, and thus capture the classical limit. 
\end{itemize}
An asymptotic-preserving method for the quantum Boltzmann equation has been proposed in~\cite{FHJ_QBE}, where a first-order
IMEX scheme combined with the standard BGK penalization idea was used. In particular, the classical Maxwellian 
was suggested in~\cite{FHJ_QBE} as an alternative to the complicated quantum Maxwellian for penalty. This replacement saves fairly amount of computational cost, but the price to pay is the loss of the strong AP property (namely, in the fluid-limit the distribution function 
should converge in one time step to its physical equilibrium state). Moreover, 
as the scheme is of IMEX type, it is hard to extend the method to very high 
order~\cite{DimarcoPareschi_HighOrderAP}. In comparison, our new schemes possess the strong AP 
property and, in principle, could achieve arbitrarily high order. As we shall
see, the presence of non classical steady states has a profound influence on 
the structure of the resulting numerical method. 

Let us finally recall that the construction of numerical methods for the full 
problem involves also discretization of the space and velocity variables. The 
latter discretization in particular is a challenging problem for the Boltzmann 
equation due to the high-dimensionality of the collision 
operator~\cite{MP_FastSpectralCollision, FHJ_QBE, HY12} and the occurrence of the Bose-Einstein 
condensation phenomenon in the degenerate quantum case~\cite{MarkPar_fast, 
MarkParBao_rev}. Here, however, we do not discuss further these issues.

The rest of the paper is organized as follows. In Section 2 we review some 
basic features of the quantum Boltzmann equation and its Euler limit. We 
emphasize in particular the differences between the classical and the quantum 
equilibrium states. Next in Section 3 we introduce the general form of the 
asymptotic-preserving exponential methods for the quantum Boltzmann equation. 
The properties of the method are then analyzed in Section 4. Several numerical 
examples are reported in Section 5 to show the AP property and the high-order 
accuracy of the schemes. We conclude the paper with some remarks in the last 
section.

\section{The quantum Boltzmann equation and its Euler limit}

The quantum Boltzmann equation was first formulated by Nordheim, Uehling and Uhlenbeck from the classical Boltzmann equation through heuristic arguments \cite{Nordheim28,UU33}. In its dimensionless form, the equation writes as:
\begin{align} \label{QBE}
\partial_tf+v\cdot \nabla_x f=\frac{1}{\varepsilon}\mathcal{Q}_q(f),\qquad t\geq 0, \quad x\in \Omega \subset \mathbb{R}^d, \quad v\in \mathbb{R}^d, \quad d=2,3,
\end{align}
where $f(t,x,v)$ is the phase space distribution function representing the (rescaled) number of particles that travel with velocity $v$ at location $x$ and time $t$. $\varepsilon$ is the so-called Knudsen number 
defined as the ratio of the mean free path over the typical length scale. It 
could vary across scales from $\varepsilon \sim \mathcal{O}(1)$ to 
$\varepsilon \ll 1$, depending on which, the system falls into the 
kinetic regime or fluid regime, respectively. The collision operator $\mathcal{Q}_q$ models the interaction between quantum particles (here and in the rest of the paper, we always use the upper sign to denote the Bose gas and the lower sign to the Fermi gas):
\begin{align}\label{eqn_Qq}
\mathcal{Q}_q(f)=\int_{\mathbb{R}^d}\int_{\mathbb{S}^{d-1}}B(v-v_*,\sigma)[f'f_*'(1\pm \theta_0f)(1\pm \theta_0f_*)-ff_*(1\pm \theta_0f')(1\pm \theta_0f_*')]\rd{\sigma}\rd{v_*},
\end{align}
where as usual, $f$, $f_*$, $f'$, and $f_*'$ are short notations for $f(t,x,v)$, $f(t,x,v_*)$, $f(t,x,v')$, and $f(t,x,v_*')$. $(v,v_*)$ and $(v',v_*')$ are the velocities before and after collision:
\begin{align}
\left\{
\begin{array}{l}
\displaystyle v'=\frac{v+v_*}{2}+\frac{|v-v_*|}{2}\sigma, \\ \\
\displaystyle v_*'=\frac{v+v_*}{2}-\frac{|v-v_*|}{2}\sigma,
\end{array}\right.
\end{align}
where $\sigma$ is the unit vector along $v'-v_*'$. The collision kernel $B$ is a nonnegative function that only depends on $|v-v_*|$ and $\cos \theta$ ($\theta$ is the angle between $\sigma$ and $v-v_*$). For variable hard sphere (VHS) particles, $B$ is independent of scattering angle:
\begin{equation}
B=C_{\gamma}|v-v_*|^{\gamma},
\end{equation} 
where $\gamma=0$ corresponds to the Maxwell molecules, and $\gamma=1$ is the hard sphere model. The parameter $\theta_0$ is some constant proportional to the Planck constant\footnote{Strictly speaking, $\theta_0=\left(\frac{2\pi \hbar}{mx_0v_0}\right)^dN$, where $m$ is the particle mass, $x_0$ and $v_0$ are the typical values of length and velocity, $N$ is the total number of particles.}:
\begin{equation}
\theta_0=C\hbar^d.
\end{equation}
It characterizes the degree of degeneracy of the system in the sense that when $\theta_0\rightarrow 0$, one recovers the collision operator for classical particles:
\begin{align}\label{eqn_Qc}
\mathcal{Q}_c(f)=\int_{\mathbb{R}^d}\int_{\mathbb{S}^{d-1}}B(v-v_*,\sigma)[f'f_*'-ff_*]\rd{\sigma}\rd{v_*}.
\end{align}

Compared with $\mathcal{Q}_c$, the quantum Boltzmann operator $\mathcal{Q}_q$ involves more nonlinearity (it is cubic rather than quadratic). This new feature brings more complexities to both theoretical and numerical studies. We are particularly interested in the fluid regime, where macroscopic equations can be derived similarly as the classical case. To this aim, we first summarize the basic properties of $\mathcal{Q}_q$.

\vspace{0.15in}
\noindent1. $\mathcal{Q}_q$ conserves mass, momentum, and energy:
\begin{align} \label{conv}
\int_{\mathbb{R}^d}\mathcal{Q}_q(f)\,\rd{v}=\int_{\mathbb{R}^d}\mathcal{Q}_q(f)v\,\rd{v}=\int_{\mathbb{R}^d}\mathcal{Q}_q(f)|v|^2\,\rd{v}=0.
\end{align}
Then if one defines the macroscopic quantities: density $\rho$, average velocity $u$, specific internal energy $e$, stress tensor $\mathbb{P}$, and heat flux $q$ as
\begin{align}
& \rho=\int_{\mathbb{R}^d}f\,\rd{v}, \quad \rho u =\int_{\mathbb{R}^d}vf\,\rd{v}, \quad \rho e=\frac{1}{2}\int_{\mathbb{R}^d}|v-u|^2f\,\rd{v}, \label{rhoue}\\
& \mathbb{P}=\int_{\mathbb{R}^{d}}(v-u)\otimes(v-u)f\,\rd{v}, \quad q=\frac{1}{2}\int_{\mathbb{R}^{d}}(v-u)|v-u|^2f\,\rd{v}, \label{Pq} 
\end{align}
the following local conservation laws can be obtained from equation (\ref{QBE}) after multiplication by $(1,v,|v|/2)^T$ and integration w.r.t. $v$:
\begin{align}
\left\{
\begin{array}{l}
\displaystyle \partial_t\rho +\nabla_x \cdot(\rho u)=0, 
\\  \\
\displaystyle \partial_t(\rho u) +\nabla_x\cdot \left(\rho u\otimes u +\mathbb{P}\right)=0, 
\\   \\
\displaystyle \partial_t\left(\rho e+\frac{1}{2}\rho u^2\right)+\nabla_x\cdot\left ( \left(\rho e+\frac{1}{2}\rho u^2\right)u+\mathbb{P}u+q\right)=0. 
\end{array}\right.
\label{local}
\end{align}

\vspace{0.15in}
\noindent2. $\mathcal{Q}_q$ satisfies the Boltzmann's H-theorem:
\begin{equation}
\int_{\mathbb{R}^d} \ln \frac{f}{1\pm \theta_0 f}\mathcal{Q}_q(f)\,\rd{v}\leq 0.
\end{equation}
Moreover, the equality holds iff $\mathcal{Q}_q(f)=0$ and iff 
$f$ reaches the local equilibrium --- the quantum Maxwellian (also called Bose-Einstein or Fermi-Dirac distribution):
\begin{equation} \label{Maxq}
\mathcal{M}_q=\frac{1}{\theta_0}\frac{1}{z^{-1}e^{\frac{(v-u)^2}{2T}}\mp1}.
\end{equation}
The new macroscopic quantities $z$ and $T$ are the fugacity and temperature. They are related to $\rho$ and $e$ via
\begin{align} \label{22system}
\left\{
\begin{array}{l}
\displaystyle  \rho=\frac{(2\pi T)^{\frac{d}{2}}}{\theta_0}Q_{\frac{d}{2}}(z),
\\   \\
\displaystyle e=\frac{d}{2}T\frac{Q_{\frac{d}{2}+1}(z)}{Q_{\frac{d}{2}}(z)},
\end{array}\right.
\end{align}
where $Q_{\nu}(z)$ is the Bose-Einstein/Fermi-Dirac function of order $\nu$ \cite{Pathria}:
\begin{align}
Q_{\nu}(z)=\frac{1}{\Gamma (\nu) }\int _0^\infty  \frac{x^{\nu-1}}{z^{-1}e^x\mp 1}dx, \quad \left( \begin{array}{ll} 0<z<1 & \text{ for Bose gas} \\  0<z<\infty  &\text{ for Fermi gas}\end{array}\right),
\end{align}
and $ \displaystyle \Gamma(\nu)=\int_0^{\infty}x^{\nu-1}e^{-x}\,\rd{x}$ is the Gamma function.

\begin{figure}[htp]
\begin{center}
   \subfigure[$\theta_0=0.01$]
   {\includegraphics[width=2.8in]{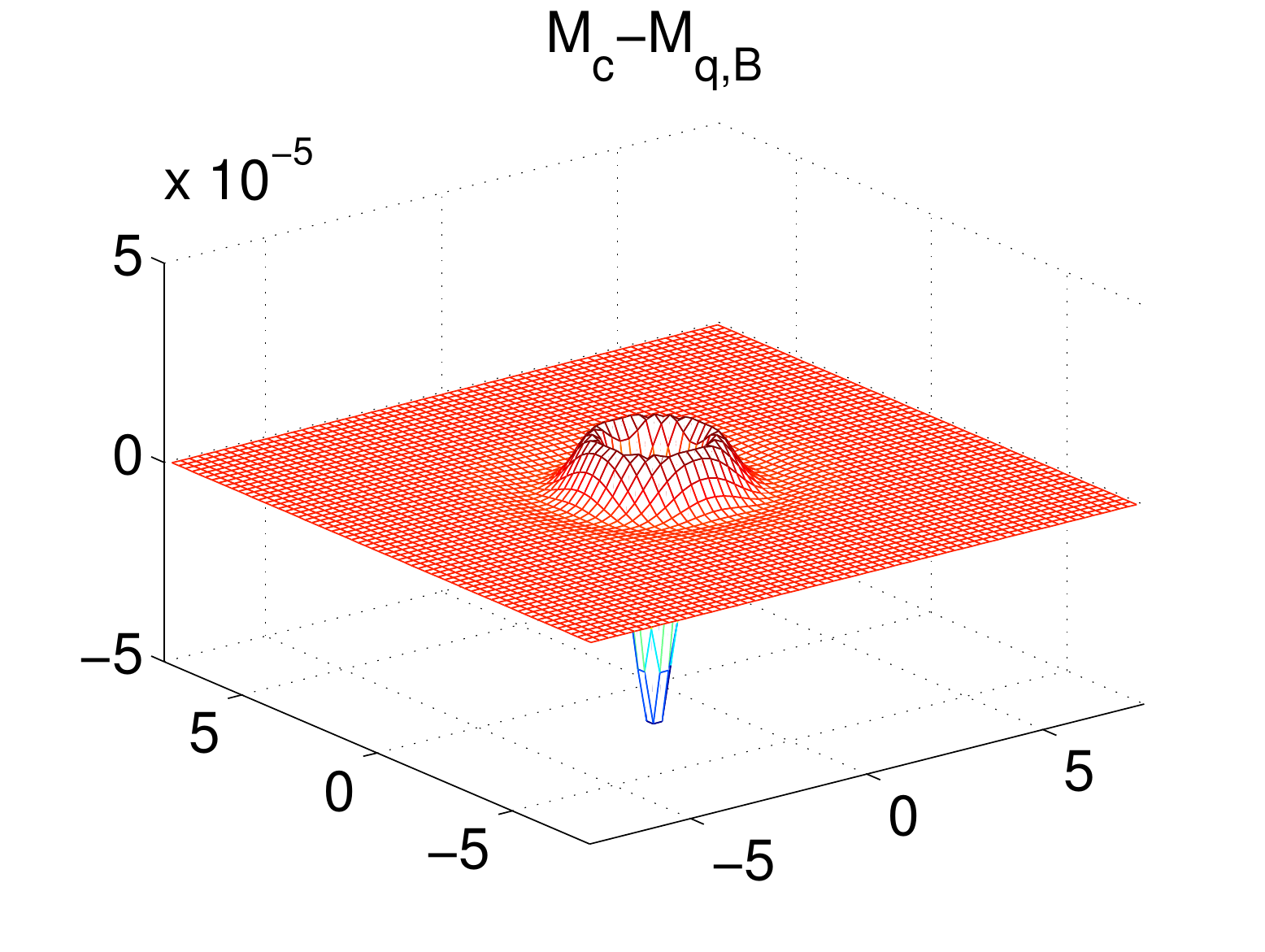}
    \includegraphics[width=2.8in]{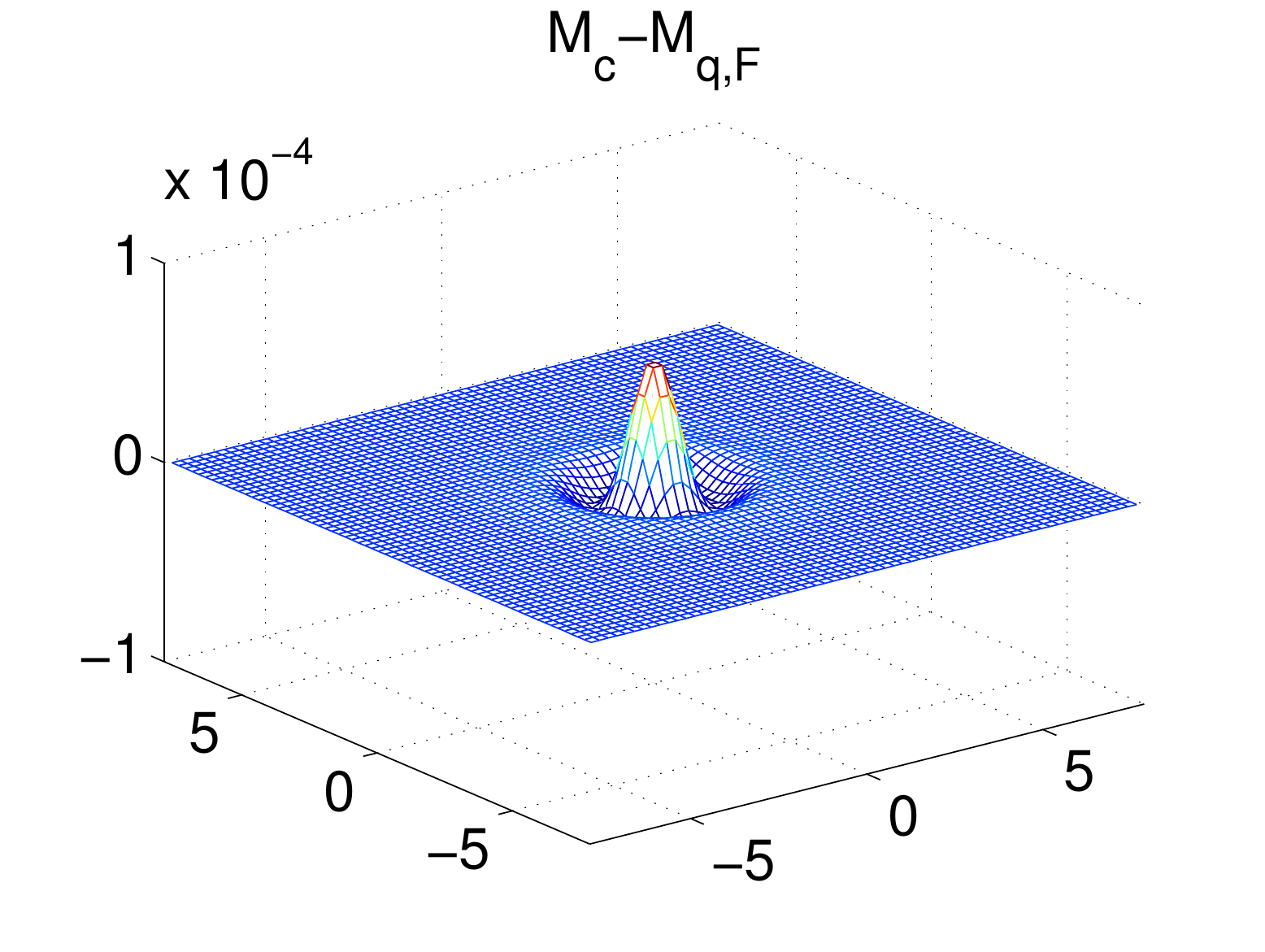}}
       \subfigure[$\theta_0=9$]
   {\includegraphics[width=2.8in]{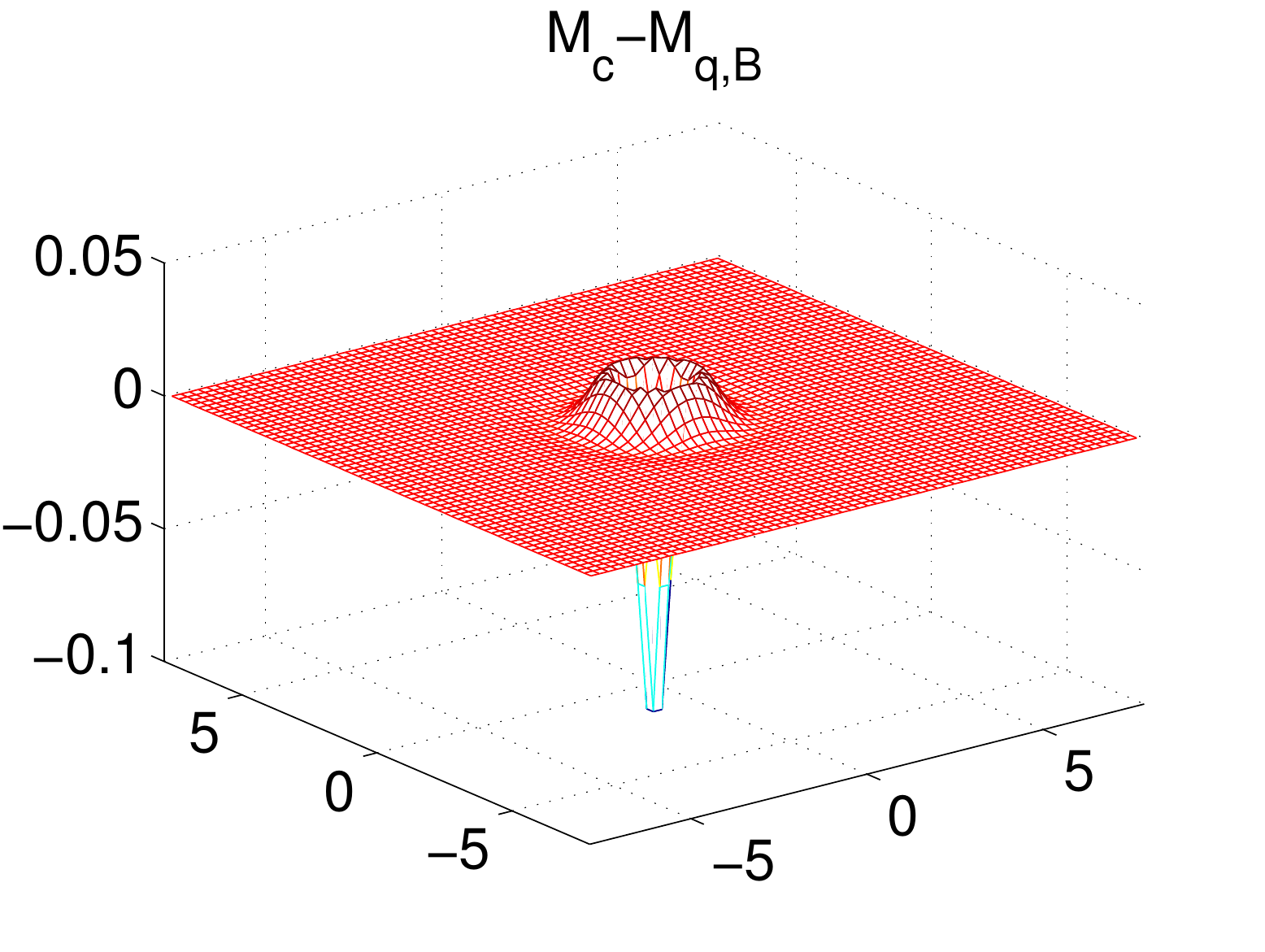}
    \includegraphics[width=2.8in]{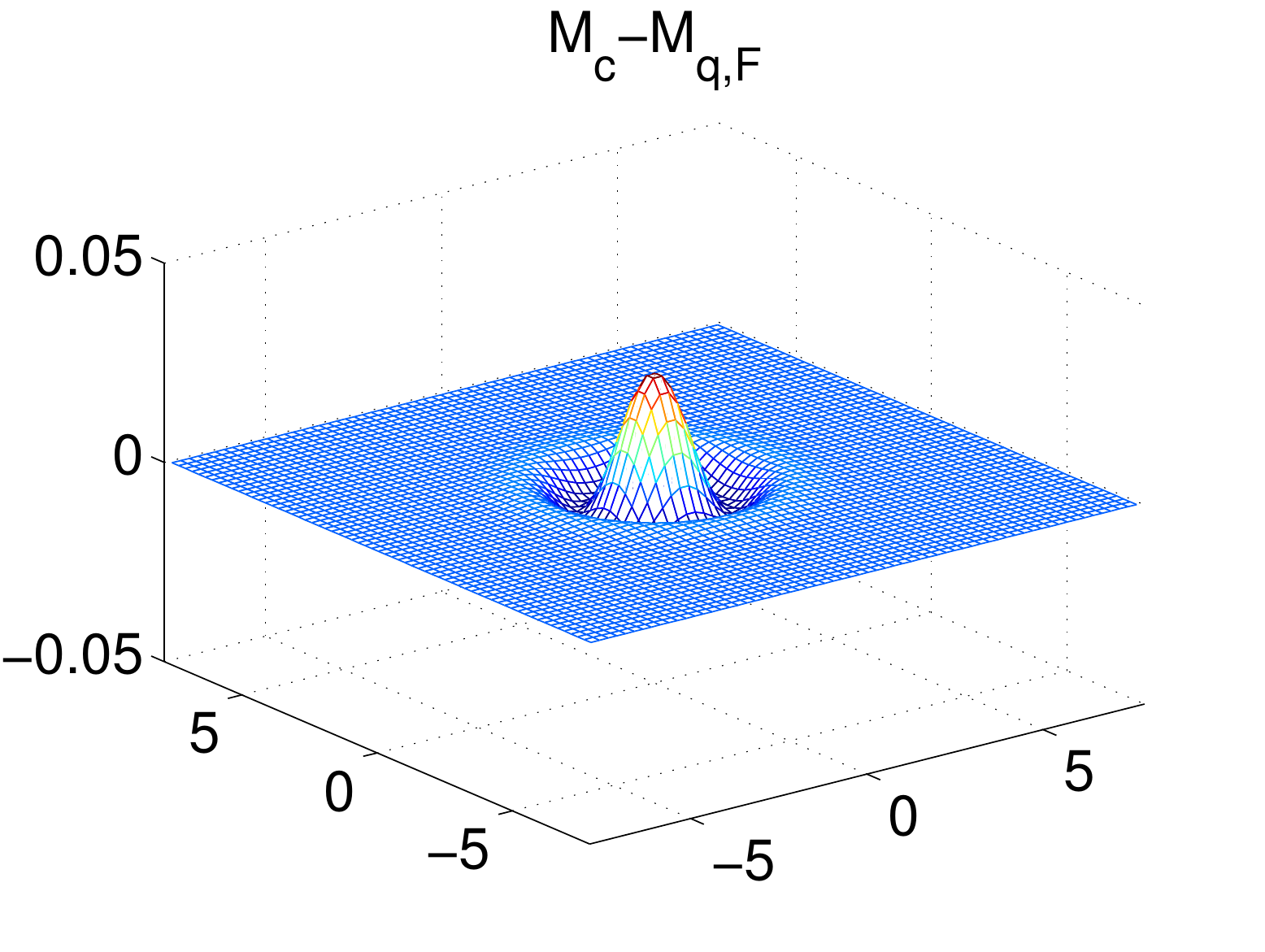}}
 \end{center}\caption{Differences between $\mathcal{M}_c$ and $\mathcal{M}_q$ in (a) the classical regime; (b) the quantum regime.
 The left column is for Bose gas and the right one is for Fermi gas. \label{fig_diff}}
\end{figure}

\begin{remark}~
\begin{itemize}
\item Compared to the classical Maxwellian:
\begin{equation}\label{Maxc}
\mathcal{M}_c=\frac{\rho}{(2\pi T)^{\frac{d}{2}}}e^{-\frac{(v-u)^2}{2T}},
\end{equation}
the quantum Maxwellian $\mathcal{M}_q$ is not a Gaussian function, and $z$, $T$ depend nonlinearly on $\rho$ and $e$, the macroscopic quantities that could be readily obtained by taking the moments of $f$. In fact, it is not difficult to see that when $z\ll 1$, $Q_{\nu}(z)$ behaves like $z$ itself. Therefore, in the system (\ref{22system}), if we keep $\rho$ and $T$ fixed, but send $\theta_0\rightarrow 0$, we get 
\begin{equation}
\frac{\rho}{(2\pi T)^{\frac{d}{2}}}\approx \frac{z}{\theta_0}, \quad e\approx\frac{d}{2}T.
\end{equation}
Since $z$ is very small, one can neglect $\mp 1$ in (\ref{Maxq}), which results in
\begin{equation}
\mathcal{M}_q\approx \frac{z}{\theta_0}e^{-\frac{(v-u)^2}{2T}}\approx \frac{\rho}{(2\pi T)^{\frac{d}{2}}}e^{-\frac{(v-u)^2}{2T}}=\mathcal{M}_c.
\end{equation}
On the other hand, if $\theta_0$ is not small, $\mathcal{M}_q$ and $\mathcal{M}_c$ will be quite different from each other. Figure~\ref{fig_diff} gives a simple illustration of the aforementioned two regimes, which we will refer to as (nearly) classical regime and quantum regime in the following discussion.
\item The physical range of interest for a Bose gas is $0<z\leq 1$, where $z=1$ corresponds to the onset of Bose-Einstein condensation (BEC). To avoid singularity, in this paper we do not consider this extreme case. We refer to~\cite{MarkPar_fast, MarkParBao_rev} for some recent results on the construction of numerical methods for the formation of BEC.
\end{itemize}
\end{remark}

\subsection{The Euler limit}

Now as the Knudsen number $\varepsilon \rightarrow 0$ in equation (\ref{QBE}), based on the discussion above, $f$ is driven to the quantum Maxwellian $\mathcal{M}_q$. Substituting $\mathcal{M}_q$ into (\ref{Pq}), we see that $\mathbb{P}=\frac{2}{d}\rho eI$ and $q=0$ ($I$ is the identity matrix). Hence the system (\ref{local}) can be closed and yields the following quantum Euler equations:
\begin{align}
\left\{
\begin{array}{l}
\displaystyle \partial_t\rho +\nabla_x \cdot(\rho u)=0, 
\\   \\
\displaystyle \partial_t(\rho u) +\nabla_x\cdot \left(\rho u\otimes u +\frac{2}{d}\rho eI\right)=0, 
\\    \\
\displaystyle \partial_t\left(\rho e+\frac{1}{2}\rho u^2\right)+\nabla_x\cdot\left ( \left(\frac{d+2}{d}\rho e+\frac{1}{2}\rho u^2\right)u\right)=0. 
\end{array}\right.
\label{Euler}
\end{align}
Obviously, written in terms of the macroscopic variables $\rho$, $u$, and $e$, this system is exactly the same as classical Euler equations. The form would be much more complicated if everything is denoted in terms of $z$, $u$, and $T$. 

\begin{remark}
By performing a Chapman-Enskog expansion to the next order, one can obtain the quantum Navier-Stokes (NS) system which differs from its classical counterpart \cite{AL97}. In particular, the viscosity and heat conductivity coefficients not only depend on $e$ but also on $\rho$. The design of a numerical scheme which is capable to capture with high accuracy the NS limit is actually under study and will be considered in future work.
\end{remark}

\section{The asymptotic-preserving (AP) exponential methods}

In this section we propose a class of high-order numerical methods for the quantum Boltzmann equation that
\begin{itemize}
\item gives accurate solution in both kinetic and fluid regimes, with time discretization not controlled by $\varepsilon$;
\item is accurate for both classical and quantum regimes with accuracy analysis independent of $\theta_0$ (or Planck constant $\hbar$).
\end{itemize}

To this end, we need to go through two steps: we 
first rewrite the equation~\eqref{QBE} in an exponential form, and then
apply the explicit Runge-Kutta methods on the newly derived equation. The 
difficulty is two fold: firstly, the equation needs to be reformulated in a way 
such that explicit Runge-Kutta, under very mild condition, automatically achieves 
AP property, and secondly, the new terms emerged in the new equation need to 
be treated with consistent schemes. We address these two difficulties in the following two subsections respectively.

\subsection{Reformulation of the equation and basic numerical methods}

In this subsection, we reformulate the equation (\ref{QBE}) in a form such that basic explicit 
Runge-Kutta methods automatically achieve asymptotic-preserving properties. The 
idea is adopted from~\cite{DP_ExpRK, LP_ExpRKinhomoBoltzmann}:
 \begin{equation}\label{eqn_transform1}
 \partial_t\left[(f-\mathcal{M}_q)e^{\frac{\mu t}{\varepsilon}}\right]=\left[\frac{\mathcal{Q}_q(f)-\mu(\mathcal{M}_q-f)}{\varepsilon}-v\cdot \nabla_x f-\partial_t\mathcal{M}_q\right]e^{\frac{\mu t}{\varepsilon}}.
 \end{equation}
This equation is derived through simple calculation, and is completely 
equivalent to the original quantum Boltzmann equation~\eqref{QBE}. However, if 
explicit methods are applied onto this equation instead of the original one, 
one could obtain the AP property.

In the equation (\ref{eqn_transform1}), $\mu$ could be any constant independent of time, and $\mathcal{M}_q$ is 
the local quantum Maxwellian function. As it shares the same moments with $f$, we 
update them through the following equation:
\begin{equation}\label{eqn_transform2}
\partial_t\int\phi{\mathcal{M}_q}\rd{v}=\partial_t\int\phi{f}\rd{v}=-\int\phi(v\cdot\nabla_xf)\rd{v}, \quad \phi=(1,v,|v|/2)^T.
\end{equation}
In the derivation we also used the fact that the first $d+2$ moments of $\mathcal{Q}_q$ are 
all zeros, as mentioned in~\eqref{conv}.

We then apply the standard Runge-Kutta method to equations (\ref{eqn_transform1}) and (\ref{eqn_transform2}). The simplest example is the forward Euler scheme:
 \begin{align}\label{scheme_FE}
 \begin{cases}
& (\displaystyle f^{n+1}-\mathcal{M}_q^{n+1})e^{\lambda}=(f^n-\mathcal{M}_q^n)+ \frac{h}{\varepsilon} \left[P^n-\mu \mathcal{M}_q^n-\varepsilon v\cdot \nabla_x f^n-\varepsilon \partial_t\mathcal{M}_q^n\right],\\
 & \displaystyle \int\phi \mathcal{M}_q^{n+1}\rd{v}=\int\phi{f}^{n+1}\rd{v}=\int\phi{f}^n\rd{v}-h\int\phi(v\cdot\nabla_xf^n)\rd{v},
 \end{cases}
 \end{align}
 where $h$ is the time step, $\lambda=\mu h/\varepsilon$, and the operator $P$ is defined as $P(f)=\mathcal{Q}_q(f)+\mu f$. The more general $\kappa$-step explicit Runge-Kutta method gives:
\begin{subequations}\label{scheme_M2}
\begin{itemize}
        \item{Stage $i$ ($i=1,\dots, \kappa$):}
            \begin{equation}\label{scheme_M2StepK}\begin{cases}
               &\displaystyle (f^{(i)}-\mathcal{M}^{(i)}_q)e^{c_i\lambda}=(f^{n}-\mathcal{M}^n_q)+\sum^{i-1}_{j=1}a_{ij}\frac{h}{\varepsilon}\left[ P^{(j)}-\mu \mathcal{M}^{(j)}_q -\varepsilon v\cdot\nabla_xf^{(j)}-\varepsilon\partial_t\mathcal{M}^{(j)}_q\right]e^{c_{j}\lambda},\\
                &\displaystyle \int \phi \mathcal{M}_q^{(i)}\rd{v}=  \int \phi f^{(i)}\rd{v} =\int\phi f^n\rd{v}+\sum^{i-1}_{j=1}a_{ij}\left(-h\int\phi v\cdot\nabla_xf^{(j)}\rd{v}\right);\end{cases}
        \end{equation}
        \item{Final Stage:}
            \begin{equation}\label{scheme_M2Final}\begin{cases}
                 & \displaystyle  (f^{n+1}-\mathcal{M}^{n+1}_q)e^{\lambda}=(f^{n}-\mathcal{M}^n_q)+\sum_{i=1}^\kappa b_i\frac{h}{\varepsilon}\left[ P^{(i)}-\mu \mathcal{M}^{(i)}_q -\varepsilon v\cdot\nabla_xf^{(i)}-\varepsilon\partial_t\mathcal{M}^{(i)}_q\right]e^{c_{i}\lambda},\\
                    &\displaystyle  \int \phi \mathcal{M}_q^{n+1}\rd{v}=           \int\phi f^{n+1}\rd{v}=\int\phi f^n\rd{v}+ \sum_{i=1}^\kappa b_i\left(-h\int\phi v\cdot\nabla_xf^{(i)}\rd{v}\right),\end{cases}
            \end{equation}
\end{itemize}
\end{subequations}
where $f^{(i)}$ stands for the estimation of $f$ at time $t=t^n+c_ih$. $a_{ij}$, $b_i$, and $c_i$ are Runge-Kutta coefficients that satisfy $\sum_{j=1}^{i-1}a_{ij}=c_i$ and $\sum_{i=1}^{\kappa}b_i=1$. They are usually stored in a Butcher tableau as:
\begin{equation}
\begin{array}{c|c}
c & A \\
\hline    
& b^T
\end{array}	
\end{equation}
Clearly at each stage $i$, to evaluate $f^{(i)}$, one needs to find $\mathcal{M}_q^{(i)}$ at the new stage first, and a good approximation of $\partial_t\mathcal{M}_q^{(j)}$ at old stages (this also applies to the final stage). 

Before moving to the next step some considerations are necessary.
\begin{remark}~
\begin{itemize}
\item The above exponential approach applied to equation (\ref{eqn_transform1}) corresponds to the so-called integrating factor method \cite{HochOstermann_ExponentialIntegrator}. Here we limit our analysis to this class of schemes, however, we refer to \cite{GPT_RelaxationNonlinearKinetic, HochOstermann_ExponentialIntegrator, LP_ExpRKinhomoBoltzmann} for other possible exponential techniques that can be used to construct other types of AP exponential schemes.

\item Reformulation (\ref{eqn_transform1}) holds true for arbitrary function $\mathcal{M}_q$ which shares the same moments with $f$. We use the local Maxwellian
$\mathcal{M}_q$ because this guarantees the strong AP property, as will be proved in Section~\ref{sec_properties}. A simplifying assumption, analyzed in \cite{LP_ExpRKinhomoBoltzmann}, consists in taking $\mathcal{M}_q$ constant along the time stepping so that the term $\partial_t\mathcal{M}_q$ disappears and the scheme simplifies. This choice, although in general less accurate in intermediate regimes, permits to obtain AP schemes with better stability and monotonicity properties. We leave the analysis of this approach in the quantum 
case for future studies and refer to \cite{LP_ExpRKinhomoBoltzmann} for further details.   
\end{itemize}
\end{remark}
 
 \subsection{Computation of $\mathcal{M}_q^{(i)}$ and $\partial_t\mathcal{M}_q^{(j)}$}
 We now show how to evaluate $\mathcal{M}_q^{(i)}$ and $\partial_t\mathcal{M}_q^{(j)}$ provided $f^{(j)}$, $\mathcal{M}_q^{(j)}$ ($j<i$) are known from previous stages. 
 
\begin{itemize}
\item[---]{Computation of $\mathcal{M}^{(i)}_q$.}\\
 By definition in~\eqref{Maxq}, $\mathcal{M}_q^{(i)}$ is obtained once we have 
 $u^{(i)}$, $z^{(i)}$, and $T^{(i)}$. The second equation in~\eqref{scheme_M2StepK} 
 gives the macroscopic quantities $\rho^{(i)}$, $u^{(i)}$, and $e^{(i)}$, and 
 thus to obtain $z^{(i)}$ and $T^{(i)}$, one only needs to invert the 
 system~\eqref{22system}. Note that the $2\times{2}$ system is nonlinear. In 
 the implementation, we use the standard Newton-iteration. Details about the 
 approximation and inversion of the quantum function $Q_{\nu}(z)$ can be found 
 in \cite{HJ}.

\item[---]{Computation of $\partial_t \mathcal{M}^{(j)}_q$.}\\
This is the key idea of the scheme. Write $\mathcal{M}_q=\mathcal{M}_q(z,T,u)$, 
it is not difficult to derive that (we drop the superscript $(j)$ for 
simplicity)
 \begin{equation} \label{Mqt}
 \partial_t \mathcal{M}_q=\mathcal{M}_q(1\pm \theta_0 \mathcal{M}_q)\left[\frac{1}{z}\partial_t z+\frac{(v-u)^2}{2T^2}\partial_t T +\frac{v-u}{T}\cdot \partial_t u\right].
 \end{equation}
While $\partial_t \rho$, $\partial_t u$, and $\partial_t e$ can be directly 
obtained from the macroscopic equations as we shall see, the computation of 
$\partial_tz$ and $\partial_tT$ is, however, not explicit. Therefore, it is 
desirable to transform the expression~\eqref{Mqt} in terms of $\partial_t \rho$, 
$\partial_t u$, and $\partial_t e$. Through the straightforward but cumbersome 
calculations in the Appendix, we end up with:
 \begin{align}\label{Mqt2}
 \partial_t \mathcal{M}_q=\mathcal{M}_q(1\pm \theta_0\mathcal{M}_q)\left[A\partial_t \rho+B\partial_te+C\cdot\partial_tu\right],
 \end{align}
 where
 \begin{align}
 A&=\frac{1}{\rho}\left(M(z)+\frac{(v-u)^2}{dT}\left(1-N(z)\right)  \right),\\
 B&=\left(\frac{(v-u)^2}{2eT}N(z)-\frac{d}{2e}M(z) \right),\\
 C&=\frac{v-u}{T},
 \end{align}
 and $M(z)$ and $N(z)$ are defined by:
 \begin{align} \label{MN}
 M(z)=\frac{Q_{\frac{d}{2}}(z)}{\left(\frac{d}{2}+1\right)Q_{\frac{d}{2}-1}(z)-\frac{d^2T}{4e}Q_{\frac{d}{2}}(z)},\quad N(z)=\frac{Q_{\frac{d}{2}-1}(z)}{\left(\frac{d}{2}+1\right)Q_{\frac{d}{2}-1}(z)-\frac{d^2T}{4e}Q_{\frac{d}{2}}(z)}.
 \end{align}

 To compute $\partial_t\rho$, $\partial_tu$, and $\partial_te$, we use 
 equation~\eqref{eqn_transform2} to transform the time derivative into spatial 
 derivative, namely:
 \begin{equation}
 \begin{cases}
\displaystyle \partial_t\rho=\partial_t\int\mathcal{M}_q\rd{v}=-\int{v}\cdot\nabla_xf\rd{v}:=-F_1,\\   
\displaystyle \partial_t\left(\rho{u}\right)=\partial_t\int{v}\mathcal{M}_q\rd{v}=-\int{v}(v\cdot\nabla_xf)\rd{v}:=-F_2,\\    
\displaystyle \partial_t\left(\rho{e}+\frac{1}{2}\rho{u}^2\right)=\partial_t\int\frac{{v}^2}{2}\mathcal{M}_q\rd{v}=-\int\frac{{v}^2}{2}({v}\cdot\nabla_xf)\rd{v}:=-F_3,
 \end{cases}
 \end{equation}
 which is
 \begin{equation}
 \begin{cases}
\displaystyle \partial_t\rho=-F_1,\\  
 \displaystyle\partial_tu=\frac{1}{\rho}\left(-F_2+F_1u\right),\\   
\displaystyle \partial_te=\frac{1}{\rho}\left(-F_3+F_1e+\frac{1}{2}F_1u^2+u\cdot(F_2-F_1u)\right).
 \end{cases}
 \end{equation}
 In this way, we could compute the time derivatives using only spatial 
 discretizations, and the scheme is automatically consistent with the 
 framework~\eqref{scheme_M2}.
 \end{itemize}

\section{Properties of the exponential AP methods}\label{sec_properties}

In this section we briefly analyze the numerical scheme. We are going to show 
that our scheme is consistent, recovers the classical Boltzmann equation in 
the classical regime, and is AP.
\begin{enumerate}
\item{\bf The classical regime:}\\
In the classical regime, $\theta_0$ (or the Planck constant) is considered as a very small 
number, and the fugacity $z\to 0$. In this regime, theoretically, the quantum 
Boltzmann equation recovers the classical Boltzmann equation, and our schemes 
should reflect this consistency.

For $z\ll 1$, as we have seen previously $Q_{\nu}(z) \approx z$, and $e \approx \frac{d}{2}T$. By definition, $M(z), N(z) \approx 1$. Plugging these relations back into
\begin{itemize}
\item{equation~\eqref{Mqt2}}: we have
\begin{equation}
A\approx \frac{1}{\rho},\quad B\approx \frac{(v-u)^2}{2eT}-\frac{d}{2e},\quad C=\frac{v-u}{T}.
\end{equation}
Therefore,
\begin{align}
 \partial_t \mathcal{M}_q \approx \mathcal{M}_q\left[ \frac{1}{\rho}\partial_t\rho + \left(\frac{d}{4e^2}(v-u)^2-\frac{d}{2e}\right)\partial_t e +\frac{d}{2e}(v-u)\cdot \partial_tu\right].
 \end{align}
 This is indeed the time evolution of the classical Maxwellian function $\partial_t\mathcal{M}_c$;
 \item{equation~\eqref{Maxq}}: as argued in Remark 2.1, $\mathcal{M}_q$ goes to $\mathcal{M}_c$;
 \item{equation~\eqref{eqn_Qq}}: formally $\mathcal{Q}_q$ also becomes the 
 classical collision operator $\mathcal{Q}_c$ as $\theta_0\to{0}$.
 \end{itemize}
 Combining these three arguments, we see that the scheme~\eqref{scheme_M2} 
 becomes the Exponential AP method developed for the classical Boltzmann 
 equation in~\cite{LP_ExpRKinhomoBoltzmann}. We successfully recovers the 
 classical regime.
 \item{\bf Consistency:}\\
 Here we assume the time step $h$ resolves $\varepsilon$. We firstly rewrite the 
 scheme~\eqref{scheme_M2StepK} as:
 \begin{equation}\label{scheme_M2StepK_expand}
 f^{(i)}=\mathcal{M}^{(i)}_q+(f^n-\mathcal{M}_q^n)e^{-c_i\lambda}+\sum_{j}a_{ij}\frac{h}{\varepsilon}\left[P^{(j)}-\mu\mathcal{M}^{(j)}_q-\varepsilon{v}\cdot\nabla_xf^{(j)}-\varepsilon\partial_t\mathcal{M}^{(j)}_q\right]e^{(c_j-c_i)\lambda}.
 \end{equation}
 As $h$ is small, we Taylor expand the exponential term, and rewrite 
 $e^{-c_i\lambda}\sim{1}-c_i\lambda$ and 
 $e^{(c_j-c_i)\lambda}\sim{1}+(c_j-c_i)\lambda$. We keep $\mathcal{O}(1)$ 
 and $\mathcal{O}(h)$ terms and neglect higher orders, the scheme becomes:
 \begin{equation}
 f^{(i)}=\Lambda_1+\Lambda_h+\mathcal{O}(h^2),
 \end{equation}
 with
 \begin{align}
&\Lambda_1=\mathcal{M}^{(i)}_q+f^n-\mathcal{M}^n_q;\\
&\Lambda_h=-c_i\lambda(f^n-\mathcal{M}_q^n)+\sum_{j}a_{ij}\frac{h}{\varepsilon}\left[P^{(j)}-\mu\mathcal{M}^{(j)}_q-\varepsilon{v}\cdot\nabla_x{f}^{(j)}-\varepsilon\partial_t\mathcal{M}^{(j)}_q\right].
 \end{align}
As $\mathcal{M}^{(i)}_q$ is the Maxwellian 
 obtained with macroscopic quantities evaluated at time $t^n+c_ih$, and thus the difference 
 $\mathcal{M}^{(i)}_q-\mathcal{M}^n_q$ is at most order $h$, therefore, one has
 \begin{equation}\label{Lambda1}
 \Lambda_1=f^n+c_ih\partial_t\mathcal{M}^n_q+\mathcal{O}(h^2).
 \end{equation}
 On the other hand, we rewrite $\Lambda_h$ as:
 \begin{align*}
 \Lambda_h=&\sum_ja_{ij}h\left(\frac{Q^{(j)}}{\varepsilon}-v\cdot\nabla_xf^{(j)}\right)\\
 &+\sum_ja_{ij}\lambda(f^{(j)}-\mathcal{M}^{(j)}_q)-c_i\lambda(f^n-\mathcal{M}^n_q)-\sum_ja_{ij}h\partial_t\mathcal{M}^{(j)}_q.
 \end{align*}
 As $f^{(j)}-f^n=\mathcal{O}(h)$, 
 $\mathcal{M}^{(j)}_q-\mathcal{M}^n_q=\mathcal{O}(h)$, and $\sum_{j}a_{ij}=c_i$, 
 we rewrite it as:
 \begin{equation}\label{Lambdah}
 \Lambda_h=\sum_ja_{ij}h\left(\frac{Q^{(j)}}{\varepsilon}-v\cdot\nabla_xf^{(j)}\right)-c_ih\partial_t\mathcal{M}^n_q+\mathcal{O}(h^2).
 \end{equation}
 Combining equation~\eqref{Lambda1} and~\eqref{Lambdah}, we have
 \begin{equation}
 f^{(i)}=f^n+h\sum_ja_{ij}\left(\frac{Q^{(j)}}{\varepsilon}-v\cdot\nabla_xf^{(j)}\right)+\mathcal{O}(h^2),
 \end{equation}
 and the consistency of the scheme is obvious. We could perform the same 
 analysis to~\eqref{scheme_M2Final} and the proof will be omitted from here.
 
 \item{\bf Asymptotic preserving:}\\
 Here we show AP property of the numerical scheme, namely, as 
 $\varepsilon\to 0$, the distribution function will automatically capture the 
 solution to the Euler equation. For simplicity, we only show proof for the case when 
 $0\leq{c}_1<c_2<\cdots<c_\kappa<1$. The argument presented here will no longer hold if any sub-stage 
 share the same time step, i.e. $c_i=c_{i+1}$ for some $i$, but we still have 
 the same conclusion. The proof for that more general case could be found 
 in~\cite{LP_ExpRKinhomoBoltzmann}.

 We still use the formula~\eqref{scheme_M2StepK_expand}. As $c_i$ 
 monotonically increases, in the zero limit of $\varepsilon$, $\lambda\to\infty$ 
 and the second and the third terms in~\eqref{scheme_M2StepK_expand} vanish, leaving:
 $$
 f^{(i)}=\mathcal{M}^{(i)}_q+\mathcal{O}\left(\lambda e^{-c\lambda}\right)\sim\mathcal{M}^{(i)}_q,\quad i=1,\ldots,\kappa,
 $$
with $c = \min_i{|c_{i+1}-c_i|}>0$. We take the moment of both sides, and combine it with the second equation in 
 the scheme~\eqref{scheme_M2StepK}:
 \begin{equation}
 \int\phi{f}^{(i)}\rd{v}\sim\int\phi{f}^n\rd{v}-\sum_{j}a_{ij}h\int\phi{v}\cdot\nabla_x\mathcal{M}^{(j)}_q\rd{v}.
 \label{eq:RKfluid}
 \end{equation}
Similarly for the numerical solution from~\eqref{scheme_M2Final} we obtain
 \begin{equation}
 \int\phi{f}^{n+1}\rd{v}\sim\int\phi{f}^n\rd{v}-\sum_{i}b_{i}h\int\phi{v}\cdot\nabla_x\mathcal{M}^{(i)}_q\rd{v}.
 \label{eq:RKfluid2}
 \end{equation}
 This is exactly how we close the moment system and obtain the Euler equation 
 analytically, and thus we capture the Euler limit. Let us note that the 
 limiting resulting scheme (\ref{eq:RKfluid})-(\ref{eq:RKfluid2}) is nothing 
 but the underlying explicit Runge-Kutta method, used in the construction of 
 the exponential scheme, applied to the limiting Euler system. Therefore the 
 method is not only consistent but it preserves the order of accuracy in the 
 fluid limit.
 \end{enumerate}
 

 \section{Numerical Examples}
  
 In this section we present several numerical results. The examples are selected to 
 reflect the AP property and high-order accuracy of the scheme we designed 
 in Section 3. Note that both examples are performed for $x$ in 1D and $v$ in 
 2D. Exp-RK2 is referred to as RK2 in time coupled with second-order 
 Lax-Wendroff scheme with van Leer limiter in space \cite{LeVeque_ConservationLaw}. Exp-RK3 is 
 referred to as RK3 in time coupled with standard WENO3 in 
 space~\cite{Shu_WENOReview}. The Butcher tableaux of RK2 (midpoint) and RK3 
 (Heun method~\cite{HaiNorWan_ODE}) we used in computation are given as follows:
 \begin{equation}
\begin{array}{c|ccc}
0   & 0   & 0   \\
1/2 & 1/2 & 0   \\
\hline
    & 0 & 1 \\
\end{array}	
\quad \quad 
\begin{array}{c|ccc}
0   & 0   & 0   & 0 \\
1/3 & 1/3 & 0 & 0  \\
2/3 & 0 & 2/3 & 0  \\
\hline
    & 1/4 & 0 & 3/4 \\
\end{array}	
\end{equation}

For the 
 velocity discretization, we use 64 points in each direction and perform 
 the fast spectral method~\cite{HY12} for Maxwell molecule kernels. Furthermore, functions $M(z)$ and $N(z)$ 
 (in the evaluation of $\partial_t \mathcal{M}_q$) have the following simple 
 form when $d=2$:
 
\begin{align}
M(z)=\frac{Q_1(z)}{2Q_0(z)-\frac{T}{e}Q_1(z)},\quad N(z)=\frac{Q_0(z)}{2Q_0(z)-\frac{T}{e}Q_1(z)},
\end{align}
where for
 \begin{itemize}
 \item{Bose gas:}  \   \  $Q_1(z)=-\ln(1-z), \quad Q_0(z)=\frac{z}{1-z};$
\item{Fermi gas:} \    \  $Q_1(z)=\ln(1+z), \quad Q_0(z)=\frac{z}{1+z}.$
 \end{itemize}
 
 \subsection{Sod problem}
\begin{figure}[htp]
 \begin{center}
 \subfigure{\includegraphics[width=3.0in]{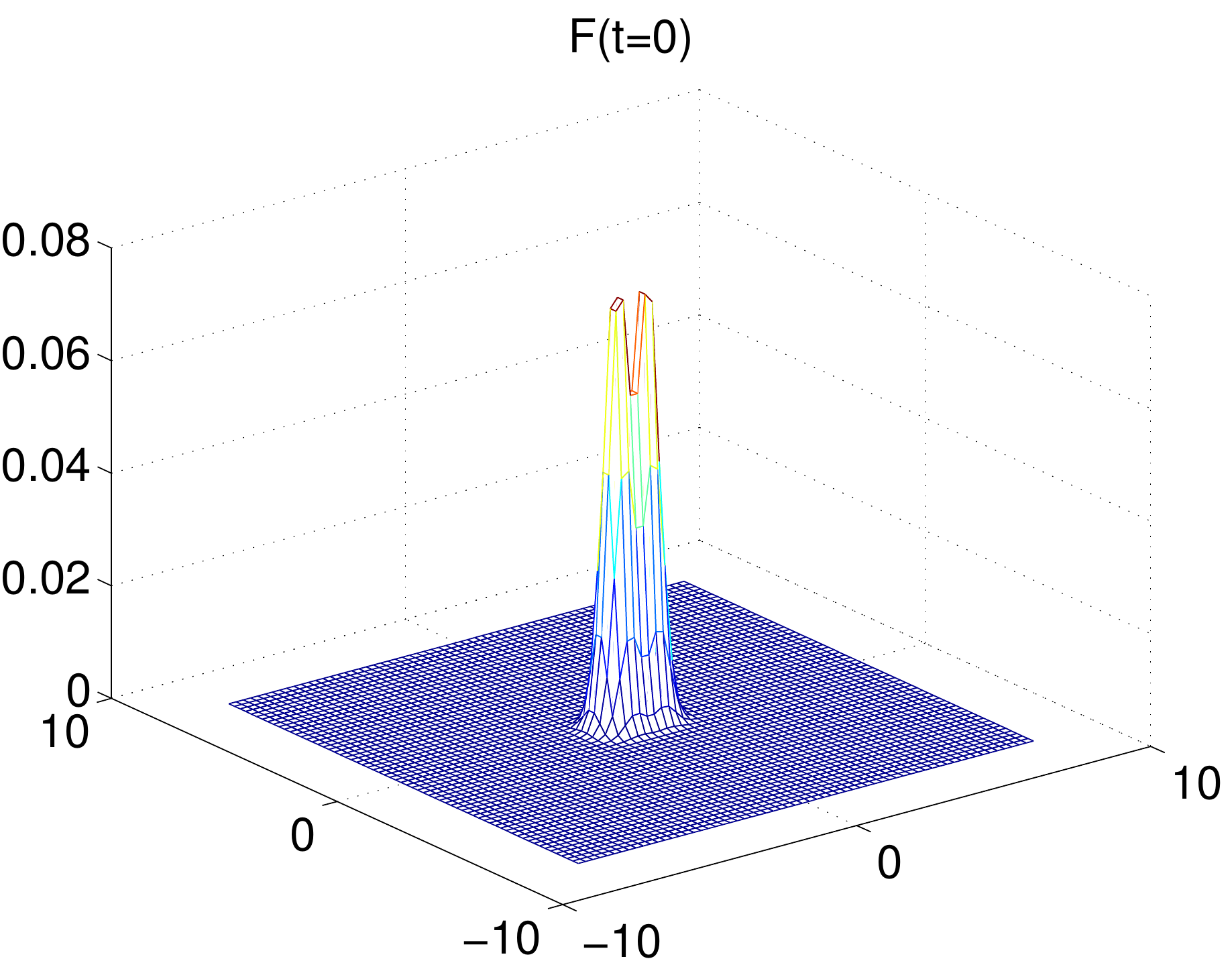}
 \includegraphics[width=3.0in]{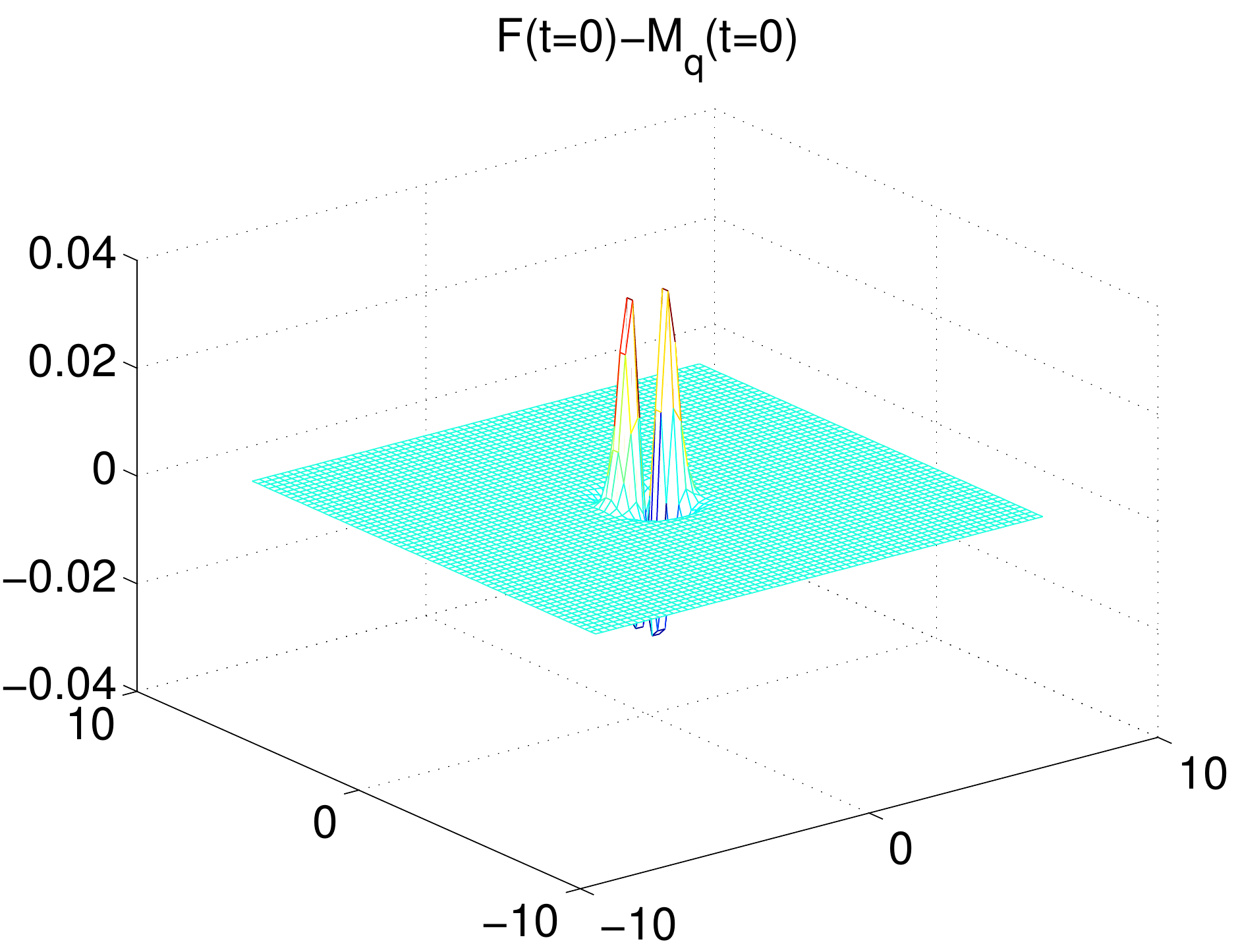}}
 \end{center}\caption{Sod problem. The initial distribution at $x=-1$. The 
 figure on left is $f(t=0)$ and right is $f-\mathcal{M}_q$.}\label{fig_iniF}
 \end{figure}

\begin{figure}[htp]
\begin{center}
   \subfigure[$\rho$]
   {\includegraphics[width=2.0in]{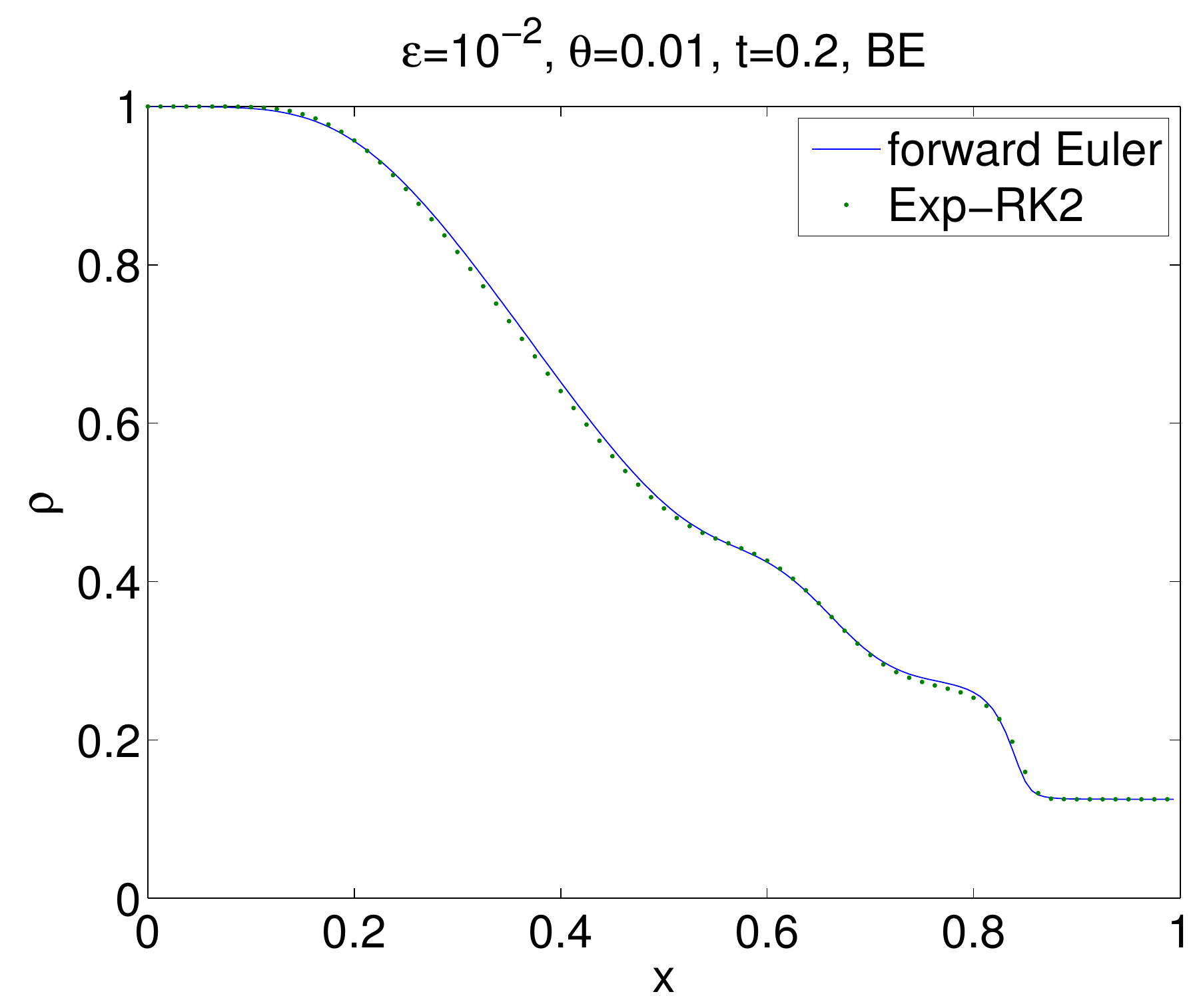}
    \includegraphics[width=2.0in]{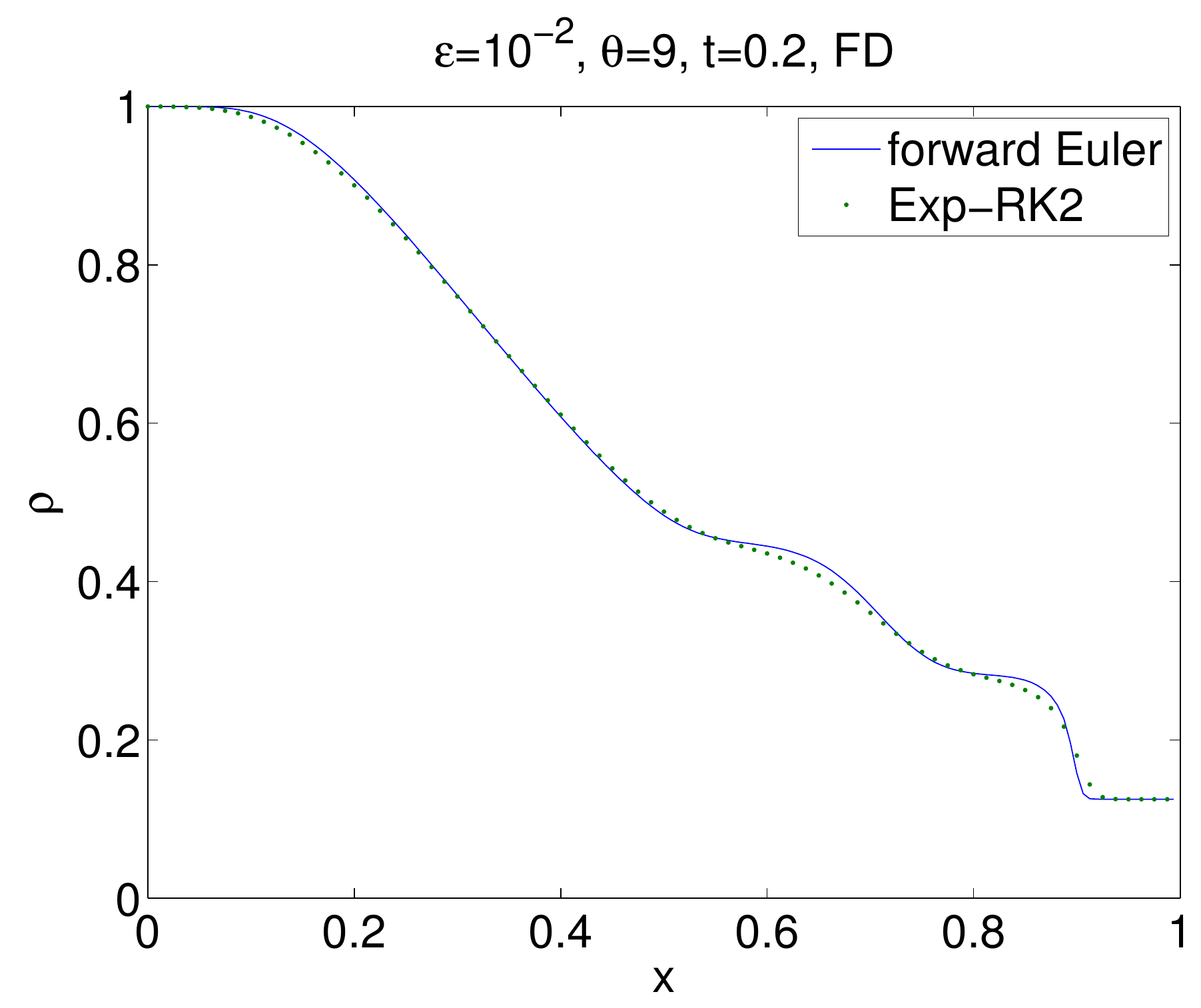}
    \includegraphics[width=2.0in]{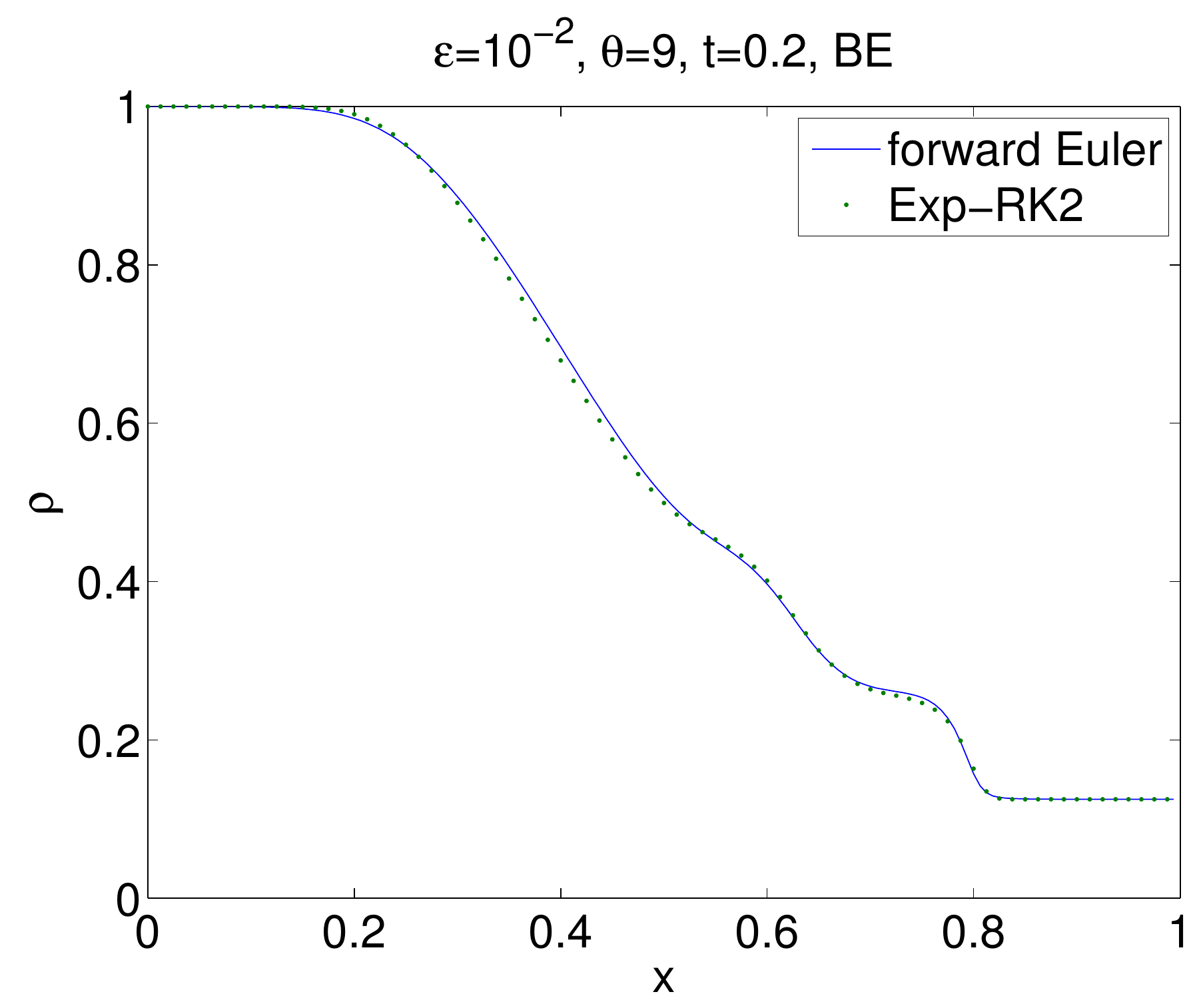}}
   \subfigure[$e$]
   {\includegraphics[width=2.0in]{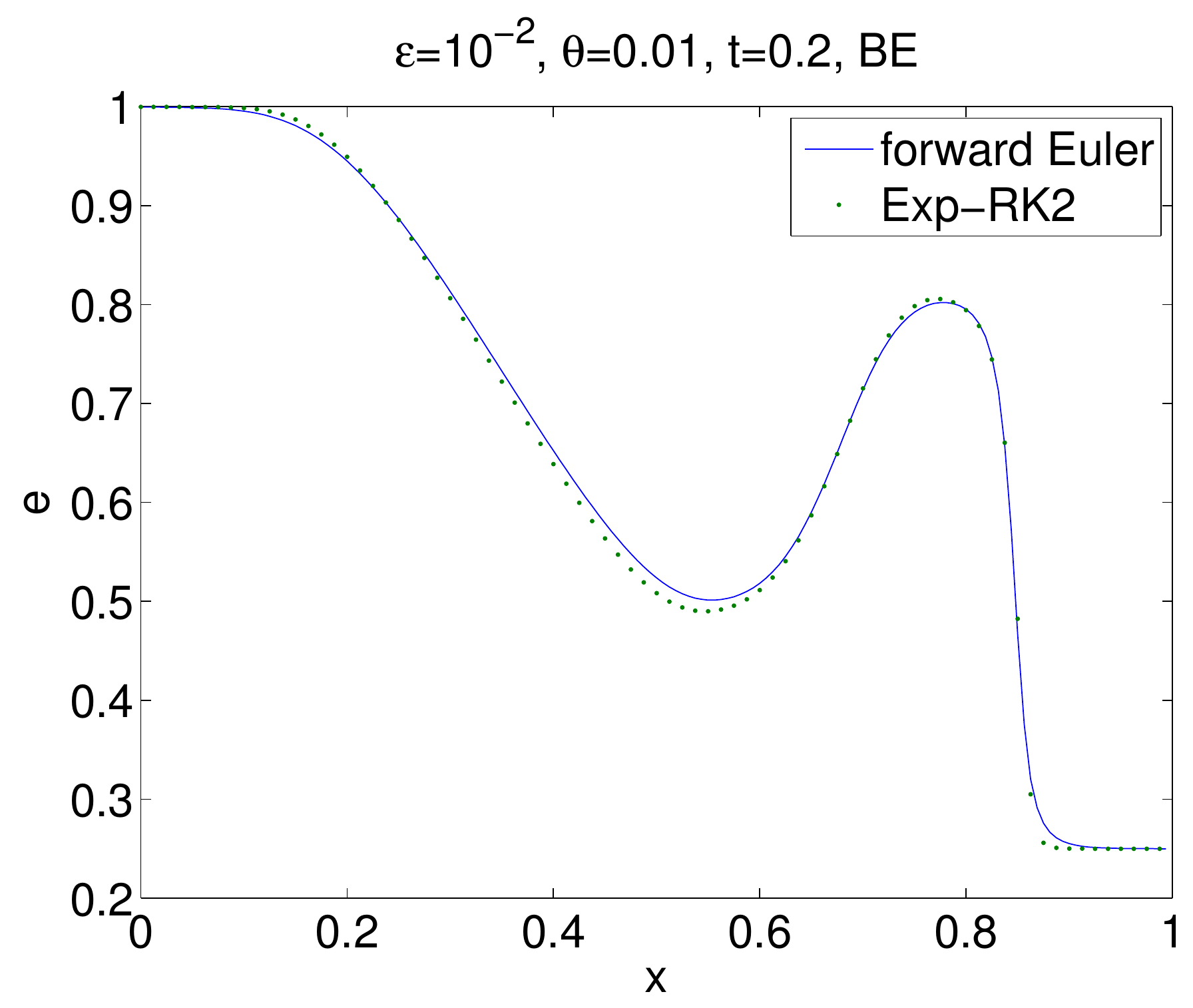}
    \includegraphics[width=2.0in]{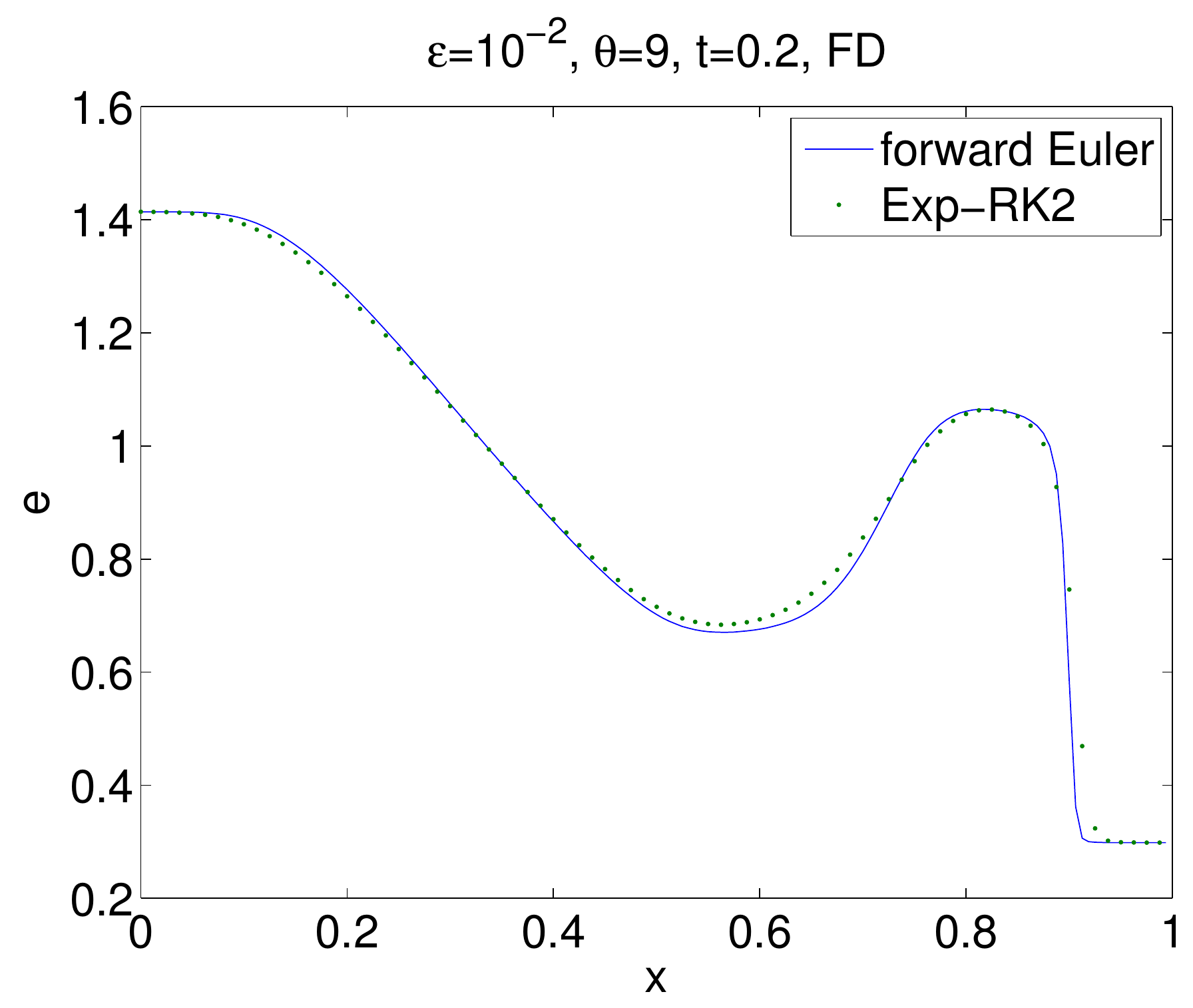}
    \includegraphics[width=2.0in]{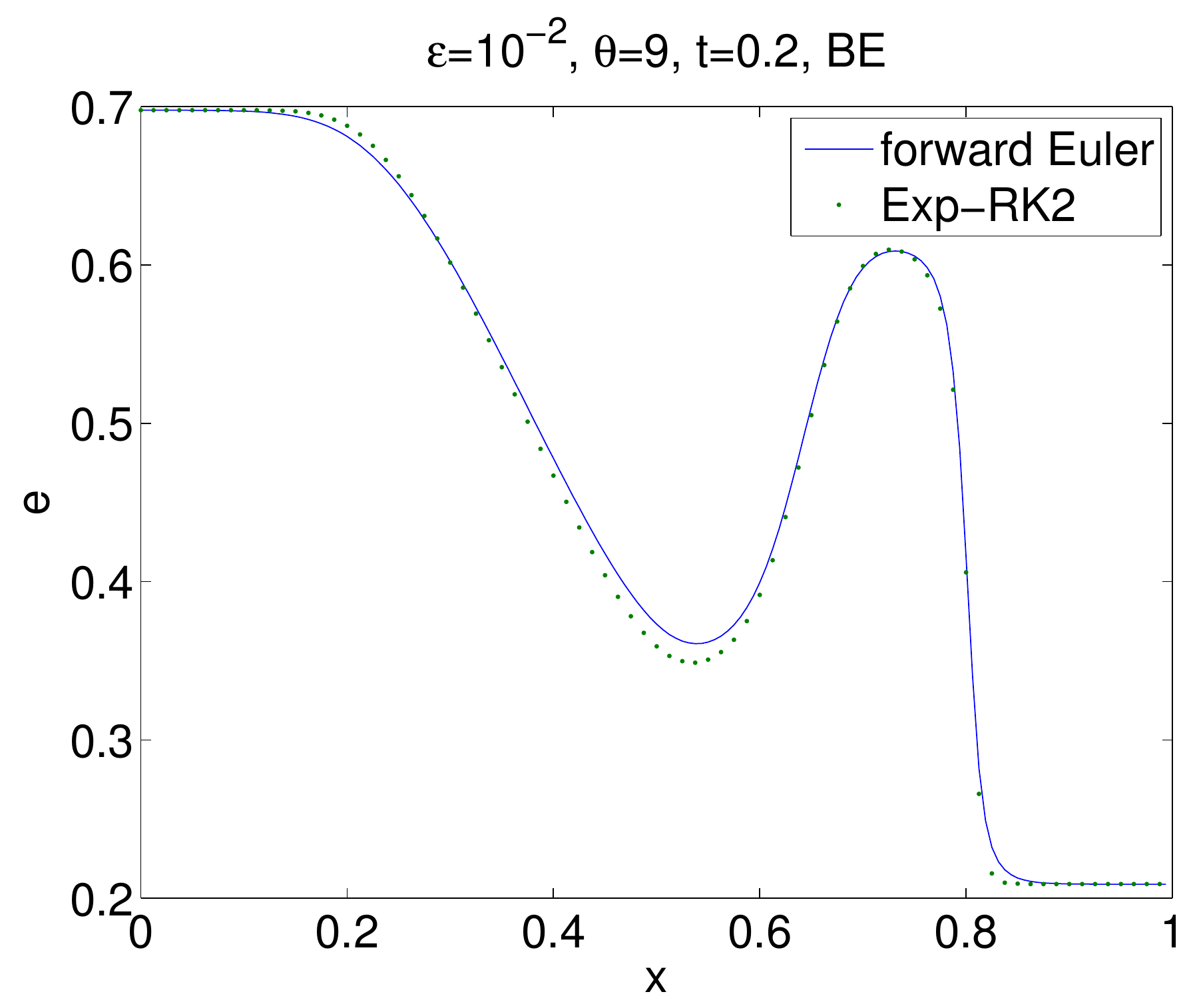}}
   \subfigure[$z$]
   {\includegraphics[width=2.0in]{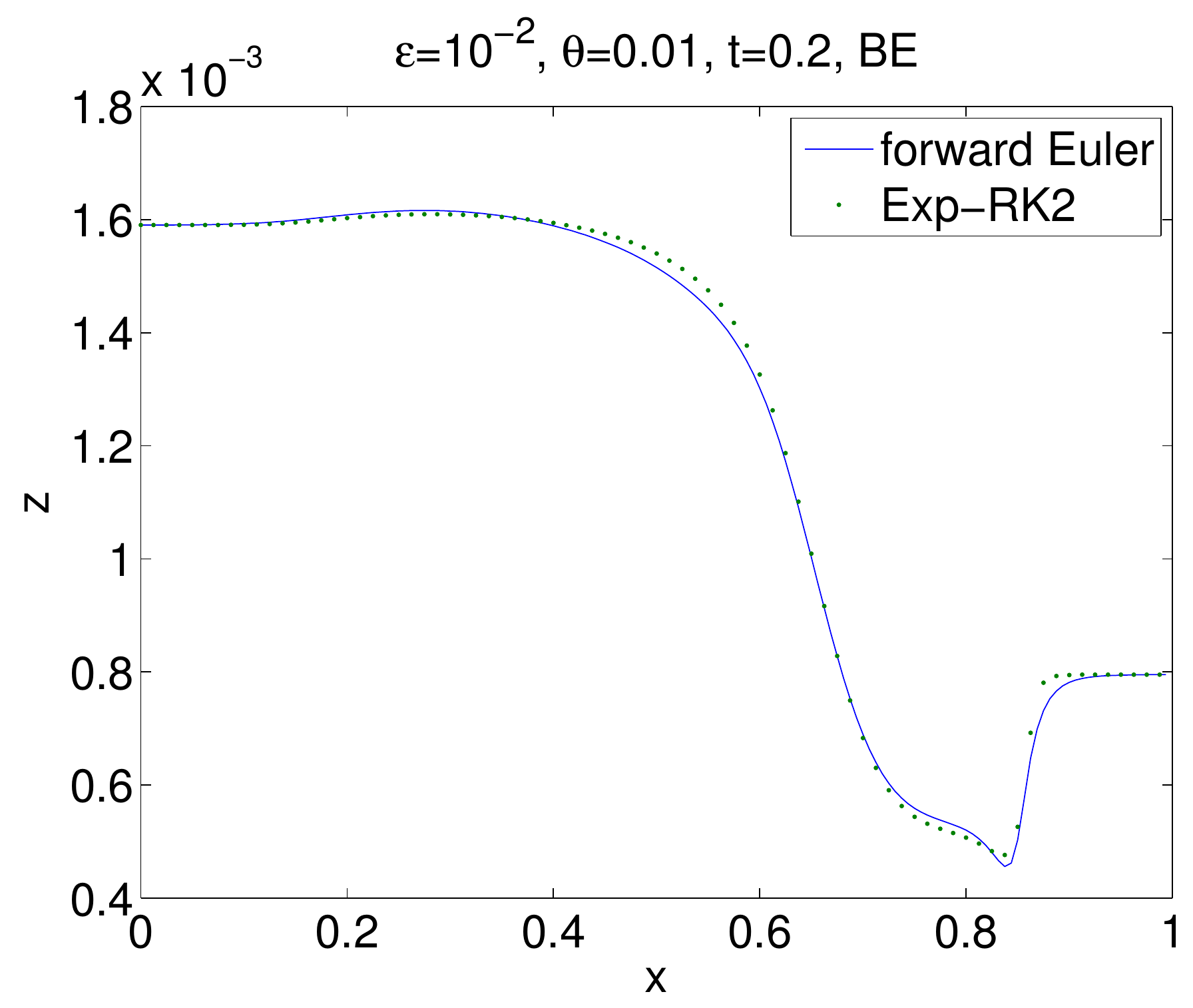}
    \includegraphics[width=2.0in]{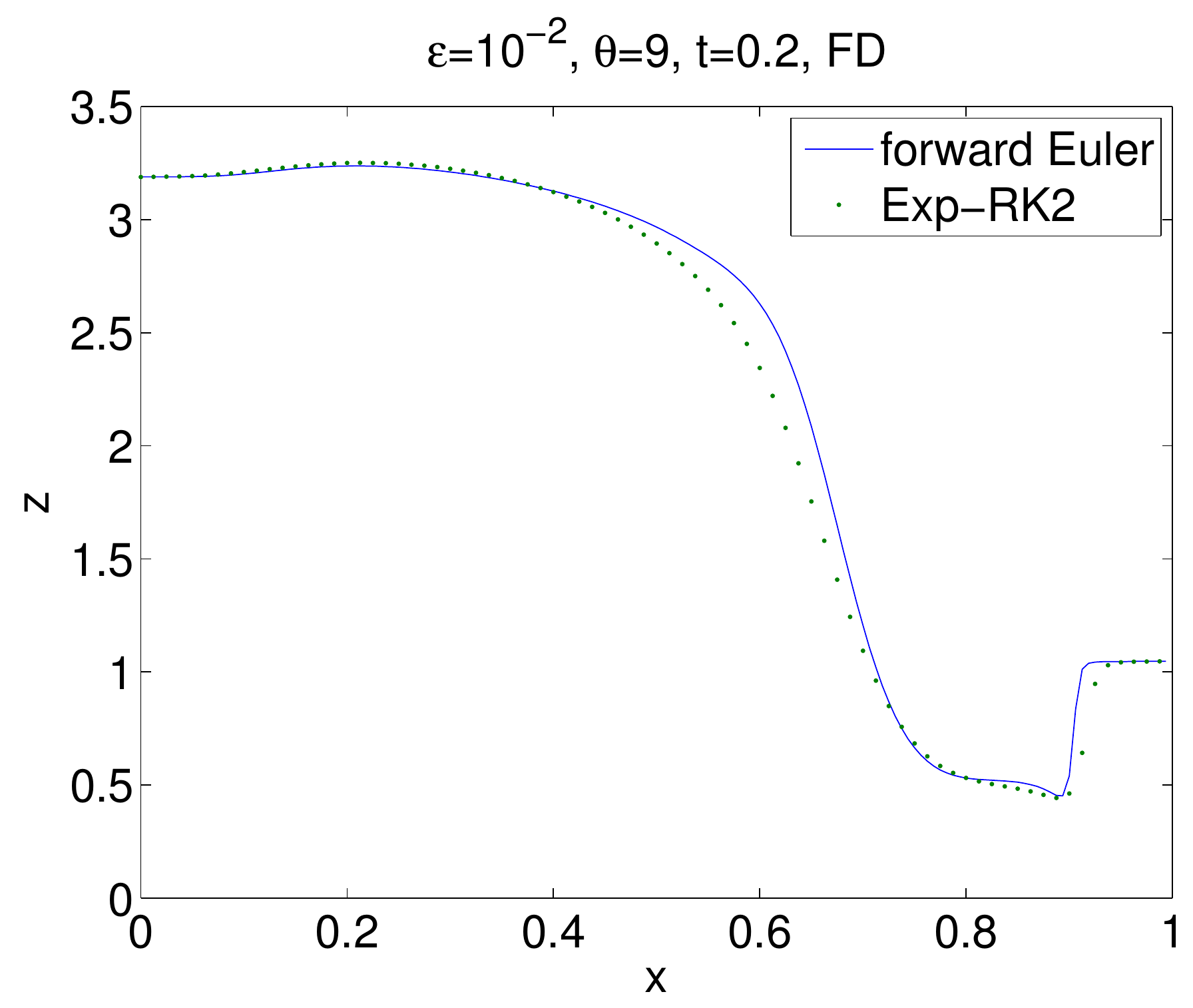}
    \includegraphics[width=2.0in]{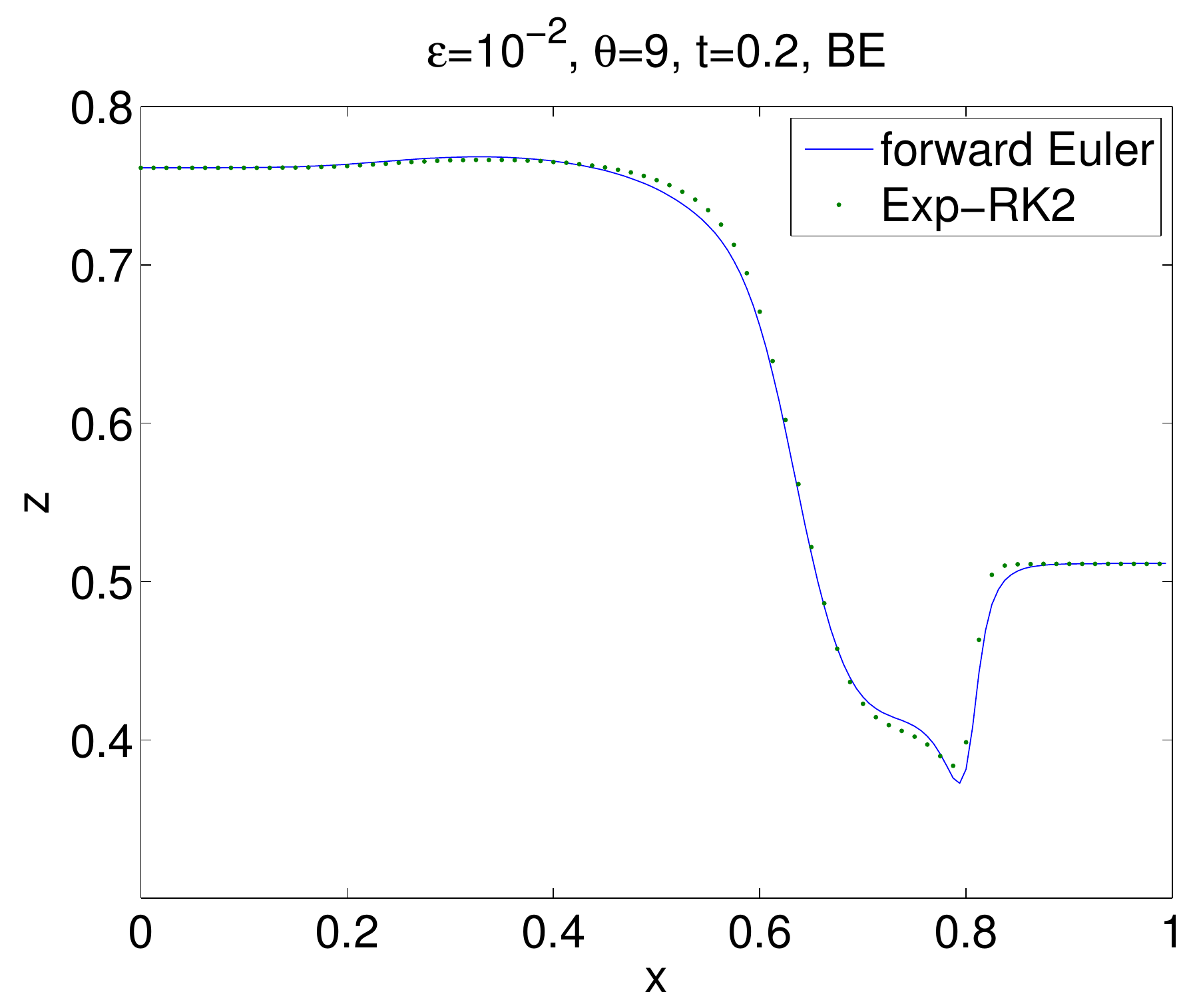}}
 \end{center}\caption{Sod problem. $\varepsilon=10^{-2}$ (kinetic regime). The three columns, from 
 the left to the right are for Bose gas in classical regime, Fermi gas in 
 quantum regime and Bose gas in quantum regime. The three rows present 
 density $\rho$, internal energy $e$ and fugacity $z$.\label{fig:Sod_ep2}}
\end{figure}

\begin{figure}[htp]
\begin{center}
   \subfigure[$\rho$]
   {\includegraphics[width=2.0in]{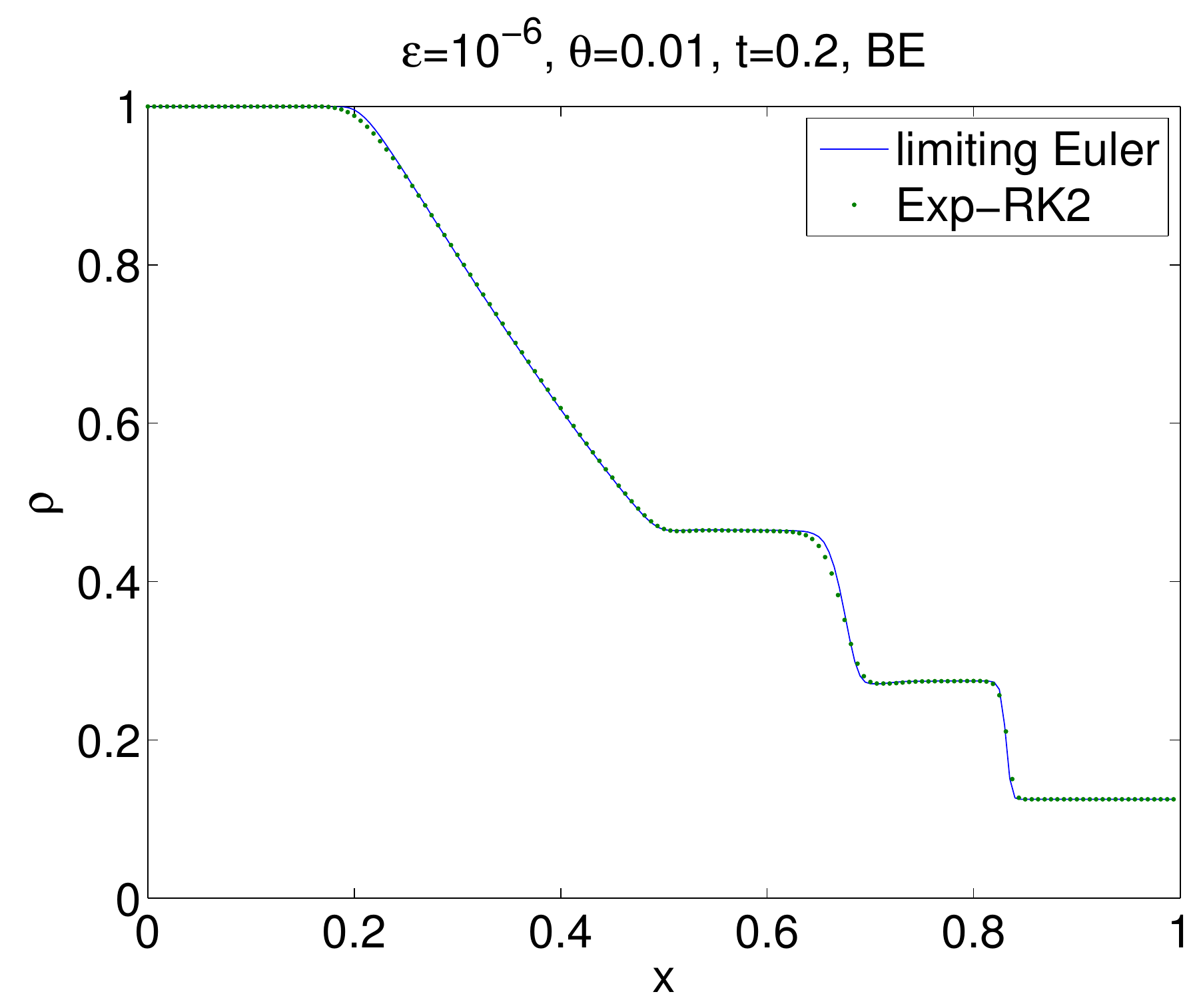}
    \includegraphics[width=2.0in]{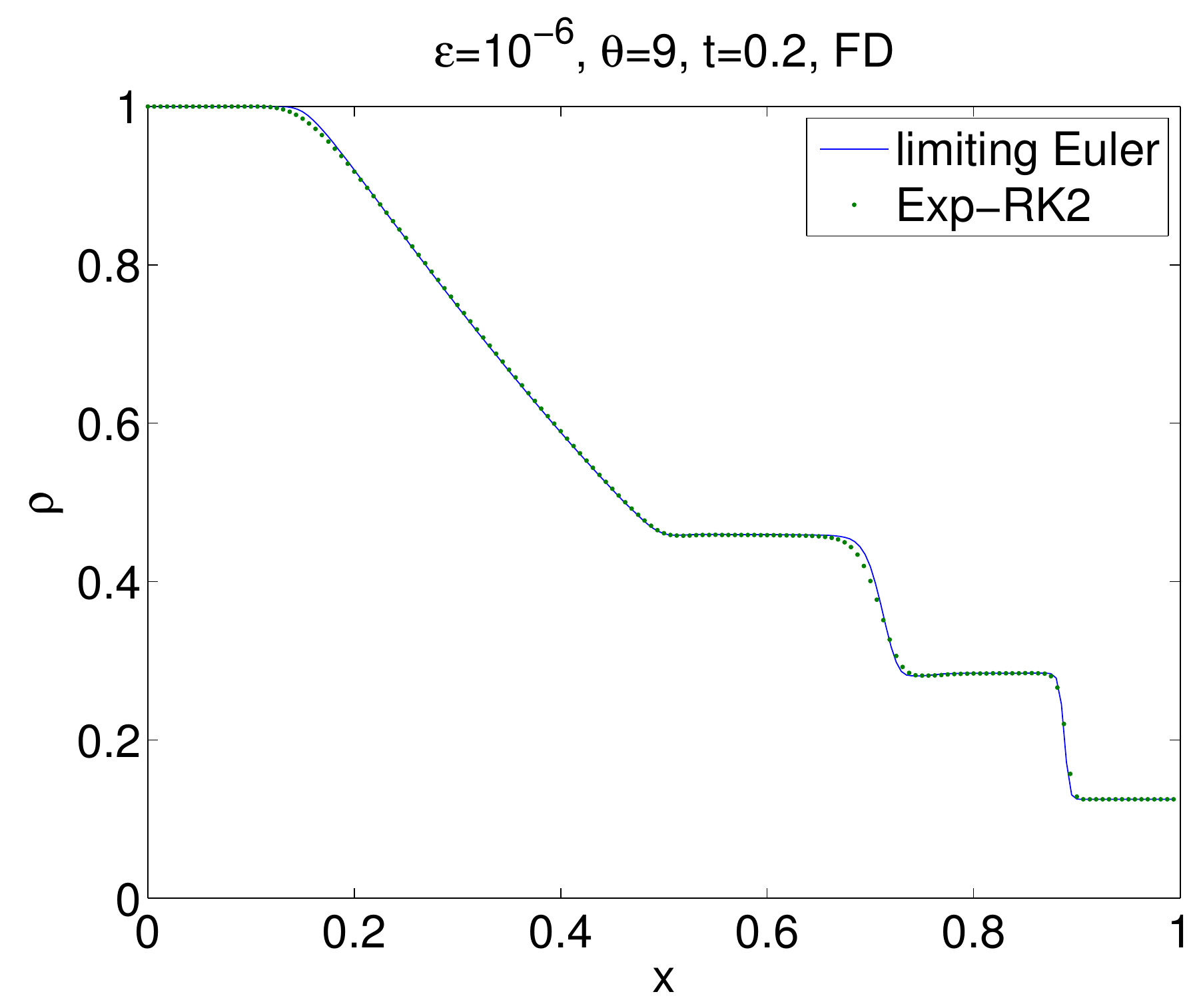}
    \includegraphics[width=2.0in]{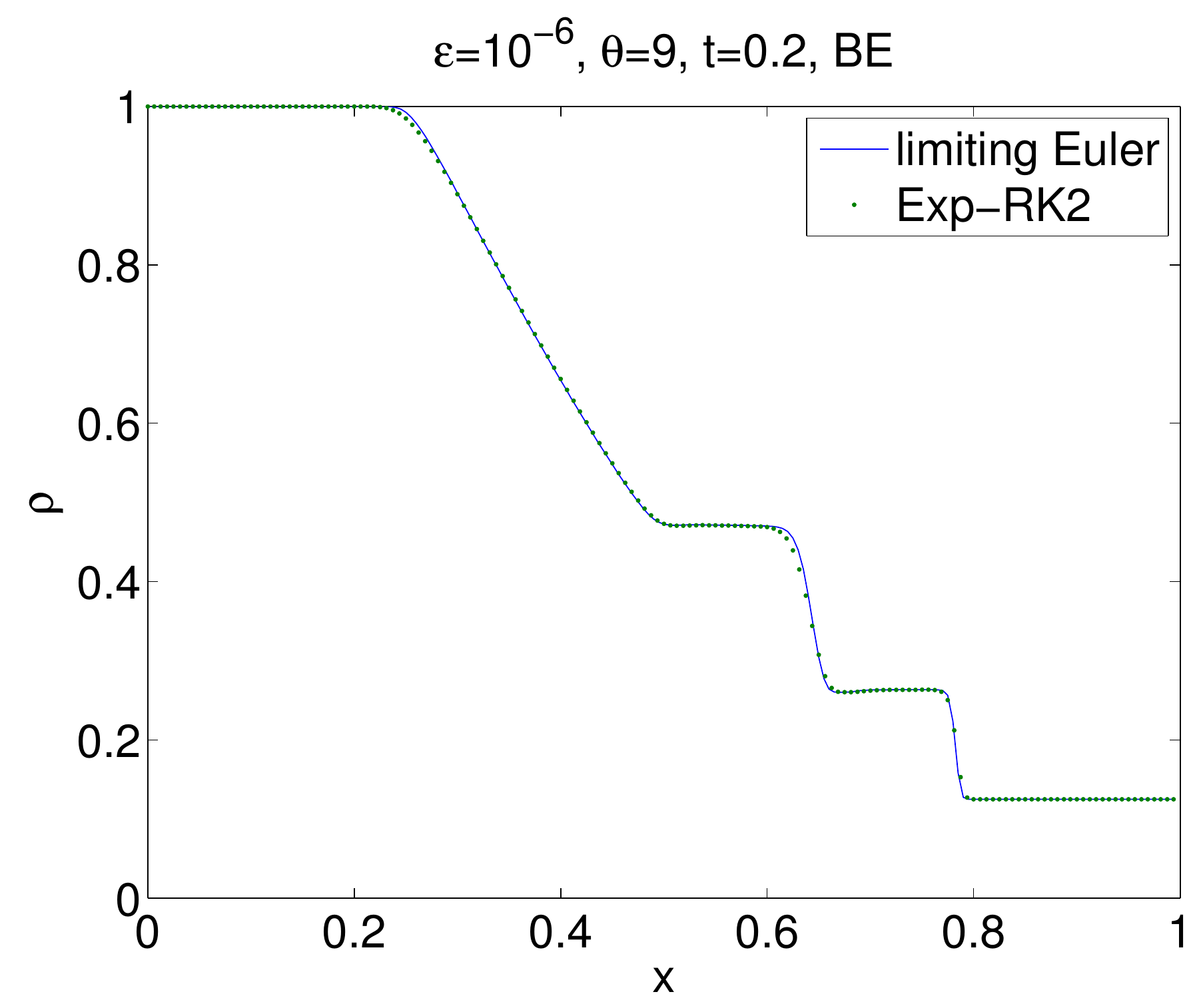}}
   \subfigure[$e$]
   {\includegraphics[width=2.0in]{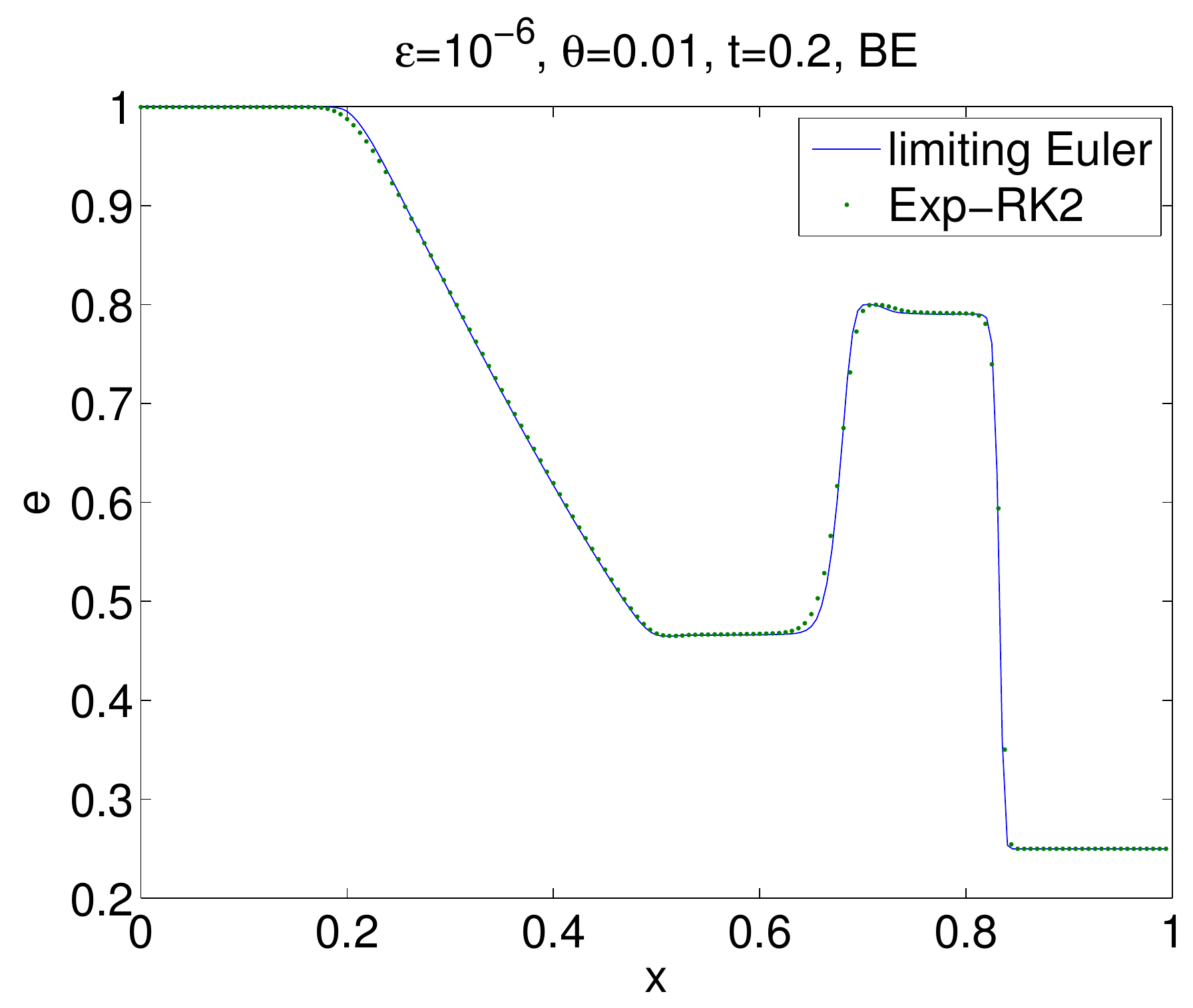}
    \includegraphics[width=2.0in]{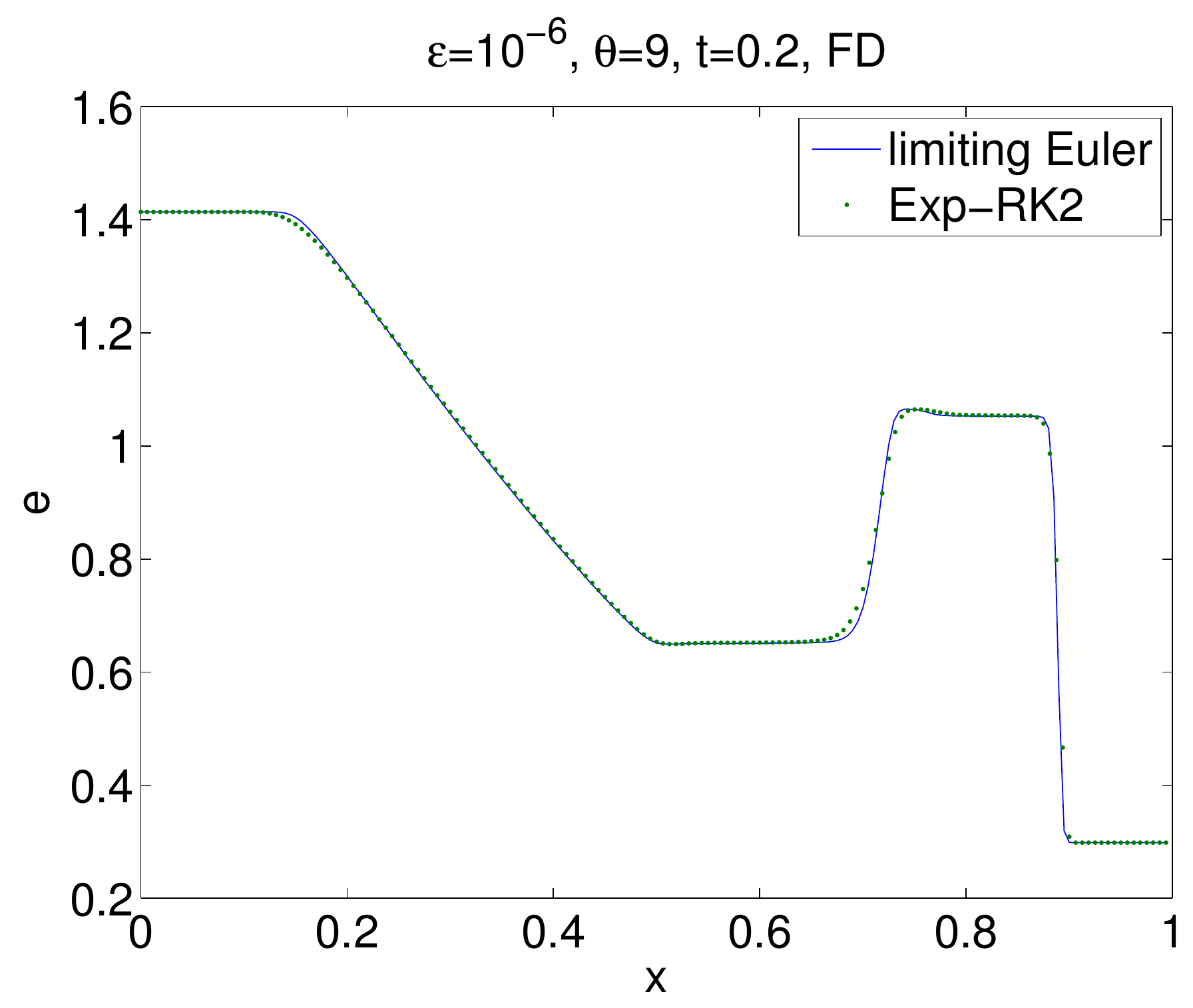}
    \includegraphics[width=2.0in]{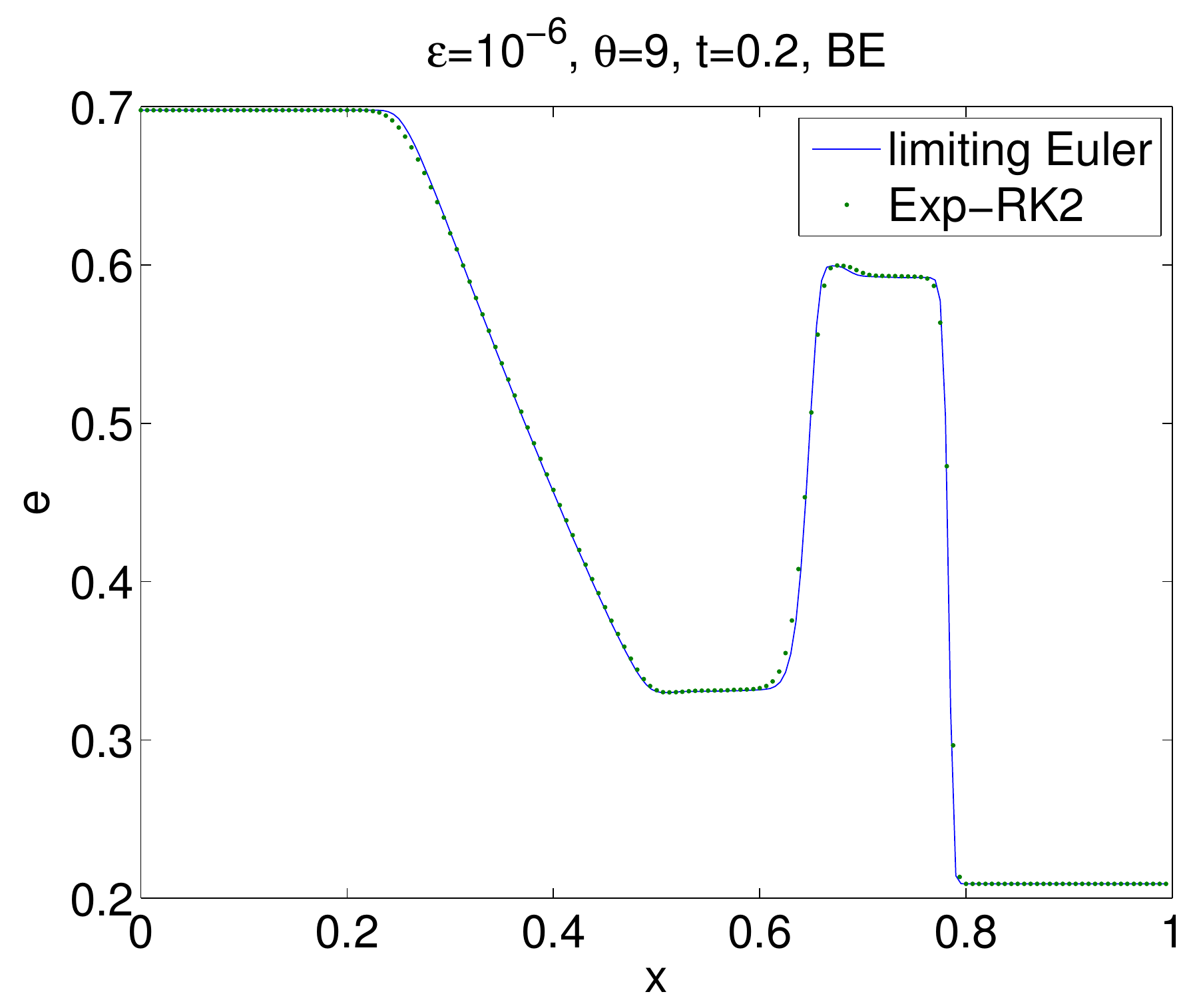}}
   \subfigure[$z$]
   {\includegraphics[width=2.0in]{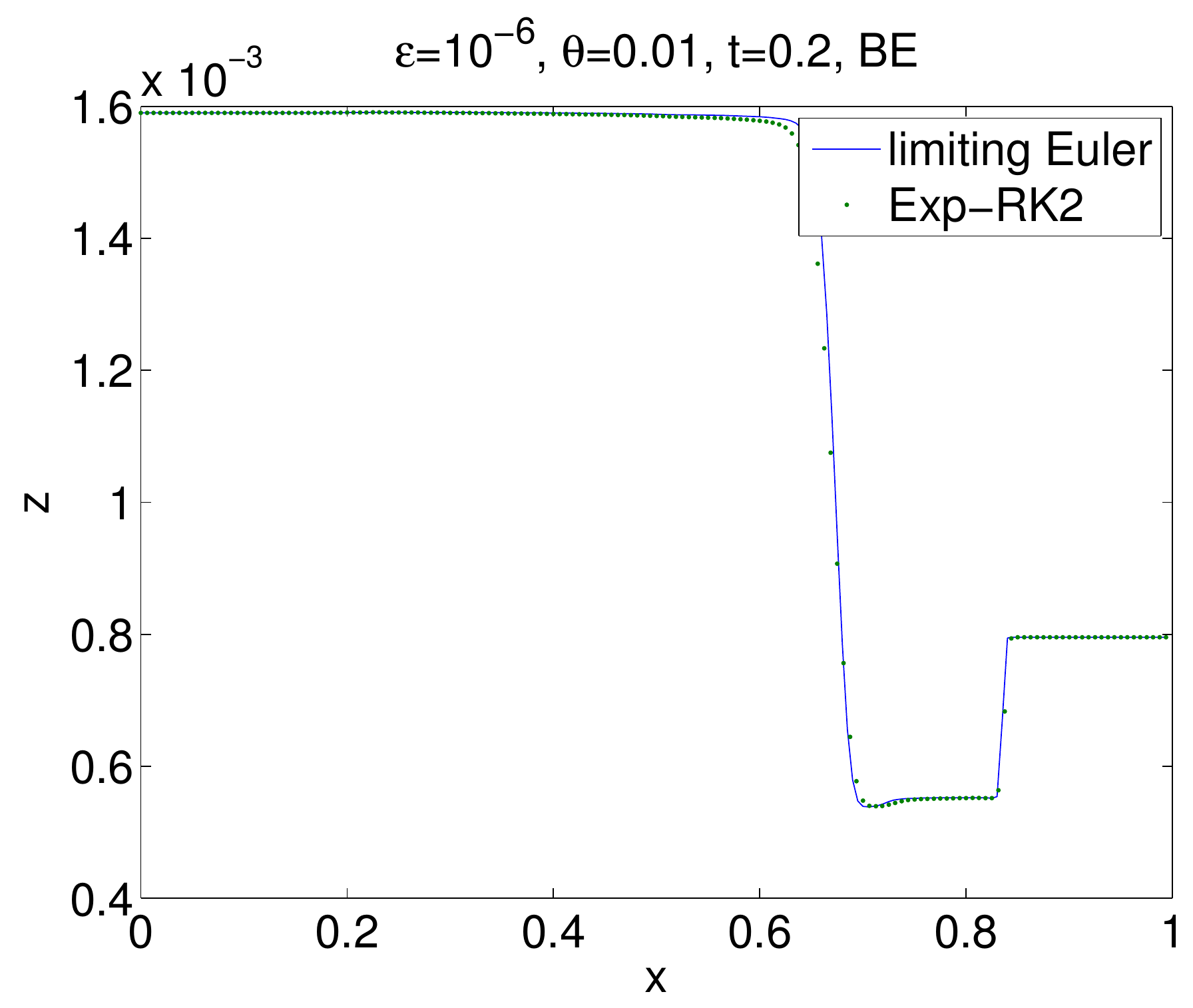}
    \includegraphics[width=2.0in]{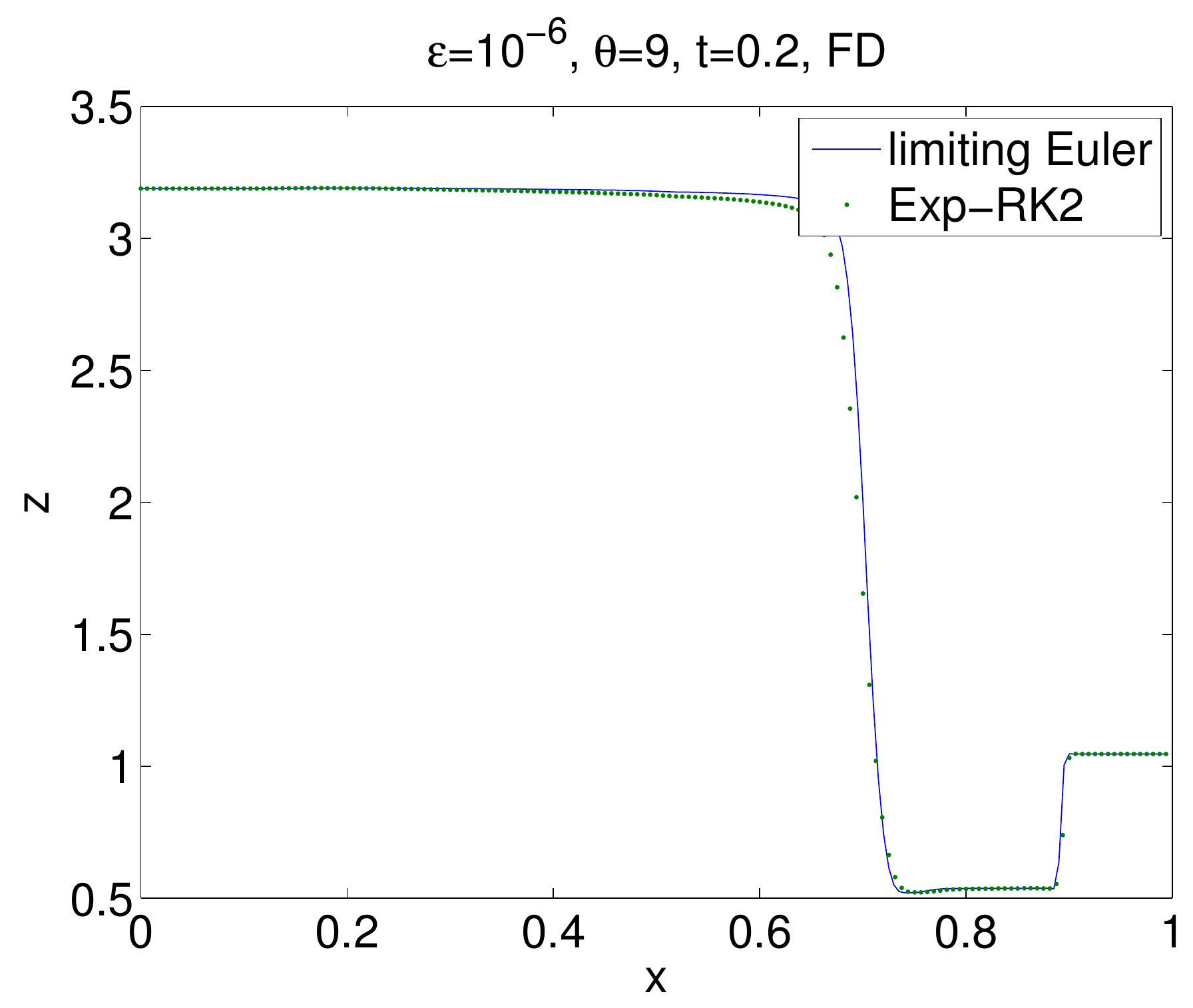}
    \includegraphics[width=2.0in]{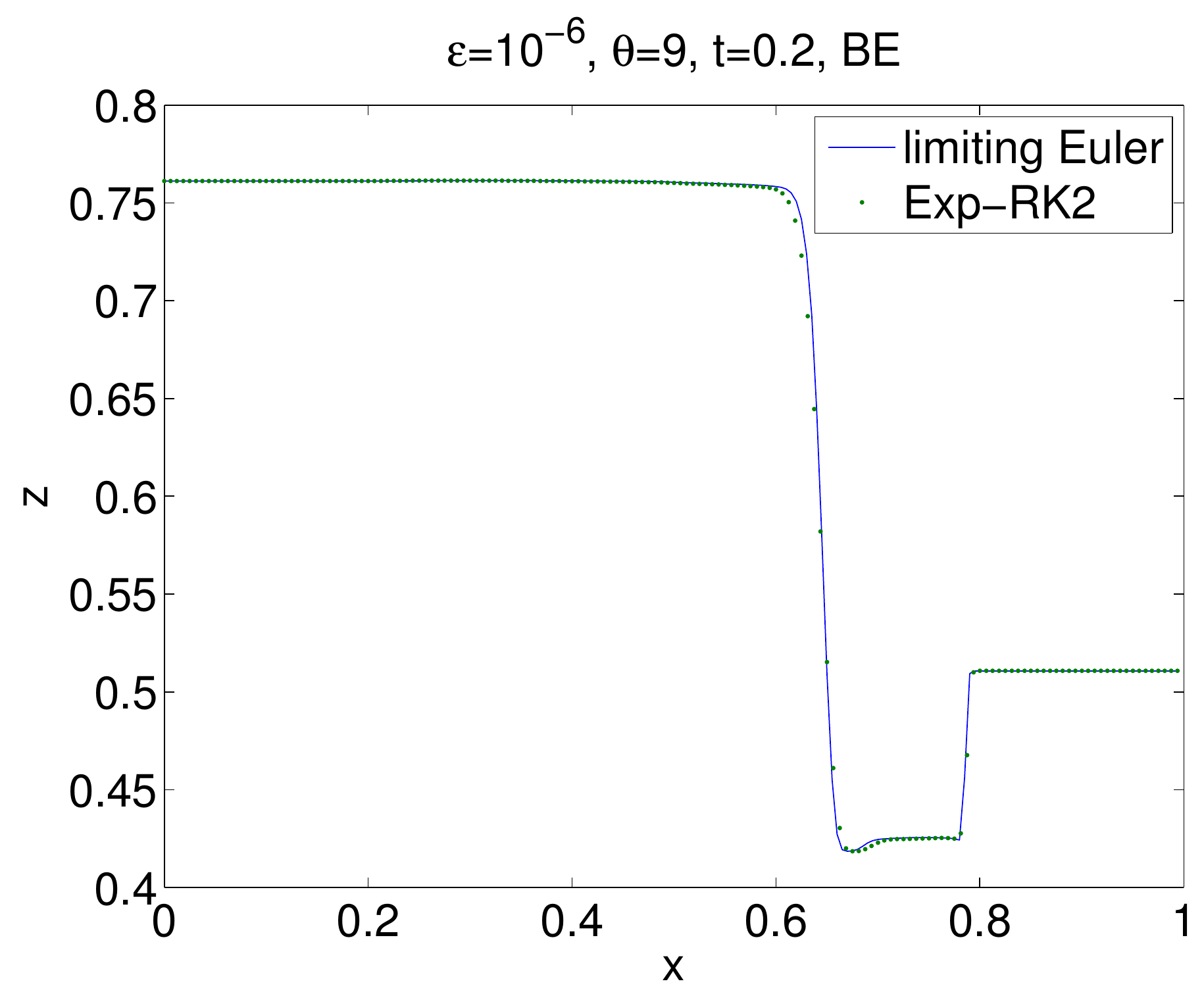}}
 \end{center}\caption{Sod problem. $\varepsilon=10^{-6}$ (fluid regime). The three columns, 
 from the left to the right are for Bose gas in classical regime, Fermi gas in 
 quantum regime and Bose gas in quantum regime. The three rows present 
 density $\rho$, internal energy $e$ and fugacity $z$.\label{fig:Sod_ep6}}
\end{figure}

 In this subsection we compute a Sod problem. In this problem, in the limiting 
 Euler regime, the solution should have a shock, a rarefaction and a contact discontinuity. The initial data for the macroscopic quantities are chosen as:
 \begin{equation}
 \begin{cases}
 \rho = 1, \quad u_x = 0, \quad u_y=0,\quad T = 1;\\
 \rho = 0.125, \quad u_x = 0, \quad u_y=0,\quad T = 0.25.
 \end{cases}
 \end{equation}
 For the microscopic quantities, we choose $f(t=0)$ to be a summation of two 
 Gaussians, as shown in Figure~\ref{fig_iniF}, and thus is far away from the quantum equilibrium.
 Figures~\ref{fig:Sod_ep2} and \ref{fig:Sod_ep6} show the numerical results using our new schemes. We consider both classical regime ($\theta_0=0.01$) and quantum regime ($\theta_0=9$) (the behaviors of Bose gas and Fermi gas in the classical regime are very close to the classical gas, thus one of them is omitted). In the case when Knudsen $\varepsilon=0.01$ (kinetic regime), the reference solutions are 
 given by directly applying the forward Euler scheme onto the original 
 Boltzmann equation with the spacial discretization $\Delta{x}=1/160$, and 
 time step $h=1/2560$, and our method uses $\Delta{x}=1/80$ and $h=1/1280$. When 
 $\varepsilon=10^{-6}$ (fluid regime), for reference data, we could not afford the fine 
 discretization any longer, and thus we directly compute the limiting Euler 
 equation. In contrast, our new scheme only uses $\Delta{x}=1/160$ and 
 $h=1/2560$, much bigger than the Knudsen number $\varepsilon$.

\begin{figure}[htp]
\begin{center}
   \includegraphics[width=2.5in]{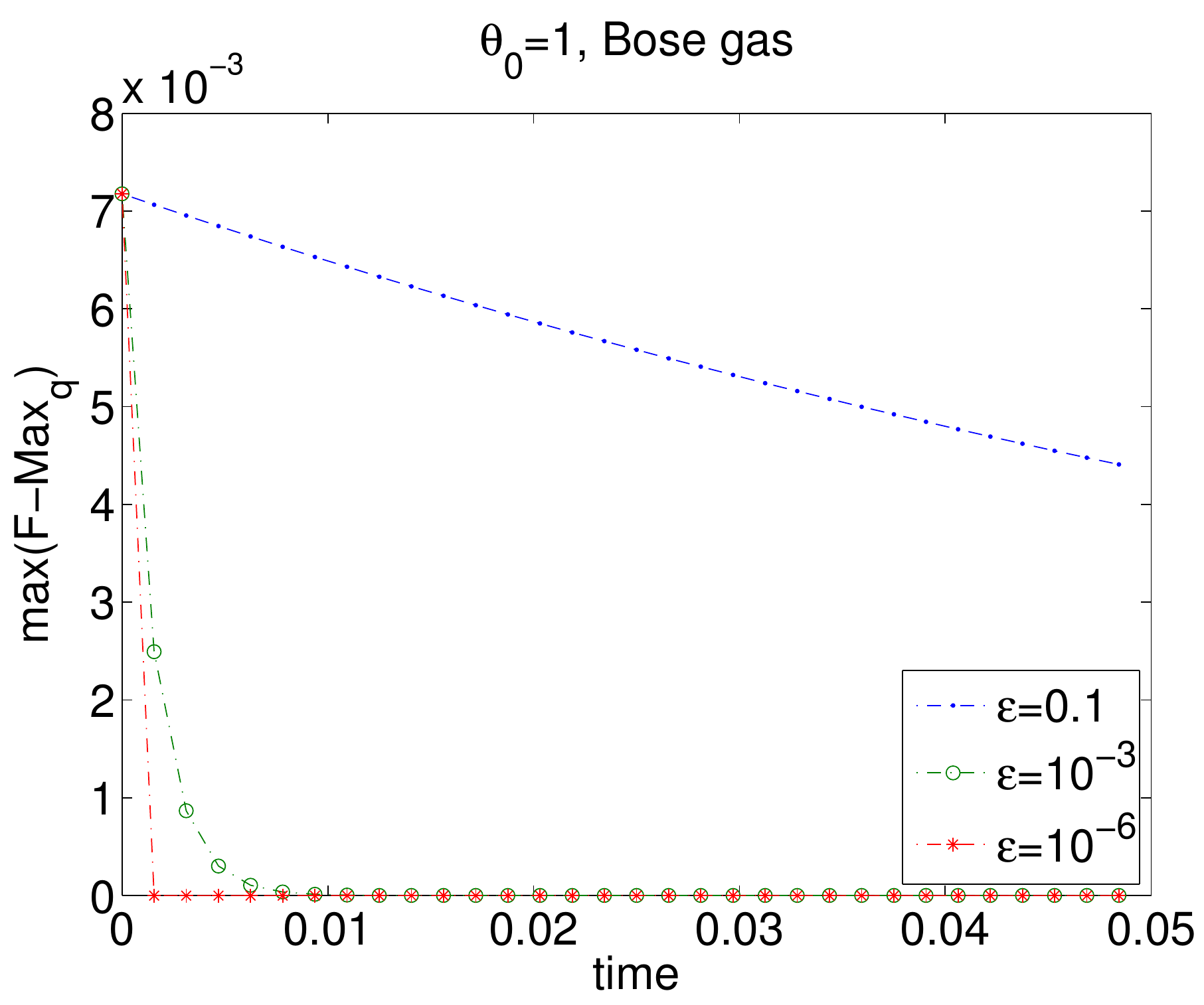}
    \includegraphics[width=2.5in]{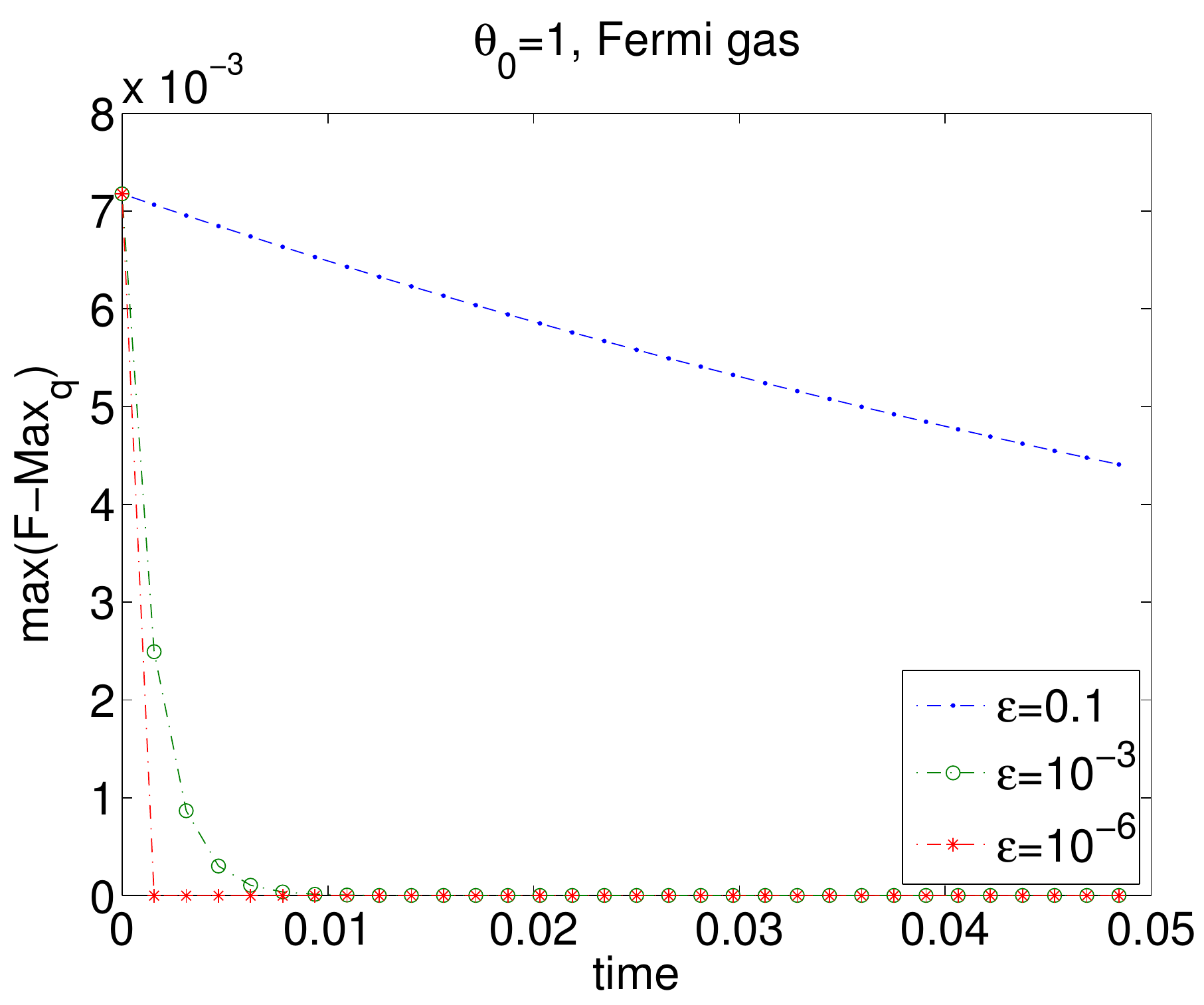}
 \end{center}\caption{Sod problem. The difference between the distribution 
 function $f$ and the Maxwellian $\mathcal{M}_q$ decays in time. The figure on 
 the left is for Bose gas and the right one is for Fermi gas. $\theta_0=1$.\label{fig_diffF}}

\end{figure}

We also measured the difference between the distribution function $f$ and the 
Maxwellian $\mathcal{M}_q$. In Figure~\ref{fig_diffF} we can clearly see that 
smaller $\varepsilon$ gives faster convergence towards the Maxwellian.

\subsection{Convergence rate test}
In the second example we show the convergence rate. We use the following smooth
initial data:
\begin{equation}
\begin{cases}
\rho = 0.3125 + 0.1875\cos{(2\pi{x})};\\
e = 0.625 + 0.375\cos{(2\pi{x})};\\
u_x = u_y = 0;
\end{cases}
\end{equation}
$h$ is chosen such that the CFL number is 
$0.5$ (independent of $\varepsilon$). Note that this is the unique stability restriction that we must impose in our numerical discretization. 

To measure the convergence rate, we check the $L_1$ error of $\rho$ and 
compute the decay rate using:
\begin{equation}
    \text{Error}_{i}=\max_{t=t^n}{\frac{\|\rho_{i}(t)-\rho_{i-1}(t)\|_1}{\|\rho_{i-1}(t)\|_1}}.
\end{equation}
Here the notation $\rho_{i}$ is $\rho$ computed on $2^i\times 20$ (with $i$ 
being a integer) grid points. 
Theoretically, if a numerical scheme is of $k$-th order, then the error should decay 
as: $\text{Error}_{i}<C\left(i\right)^{-k}$ for $h$ small enough.

\begin{figure}[htp]
\begin{center}
   \subfigure[Bose gas]
   {\includegraphics[width=3.2in]{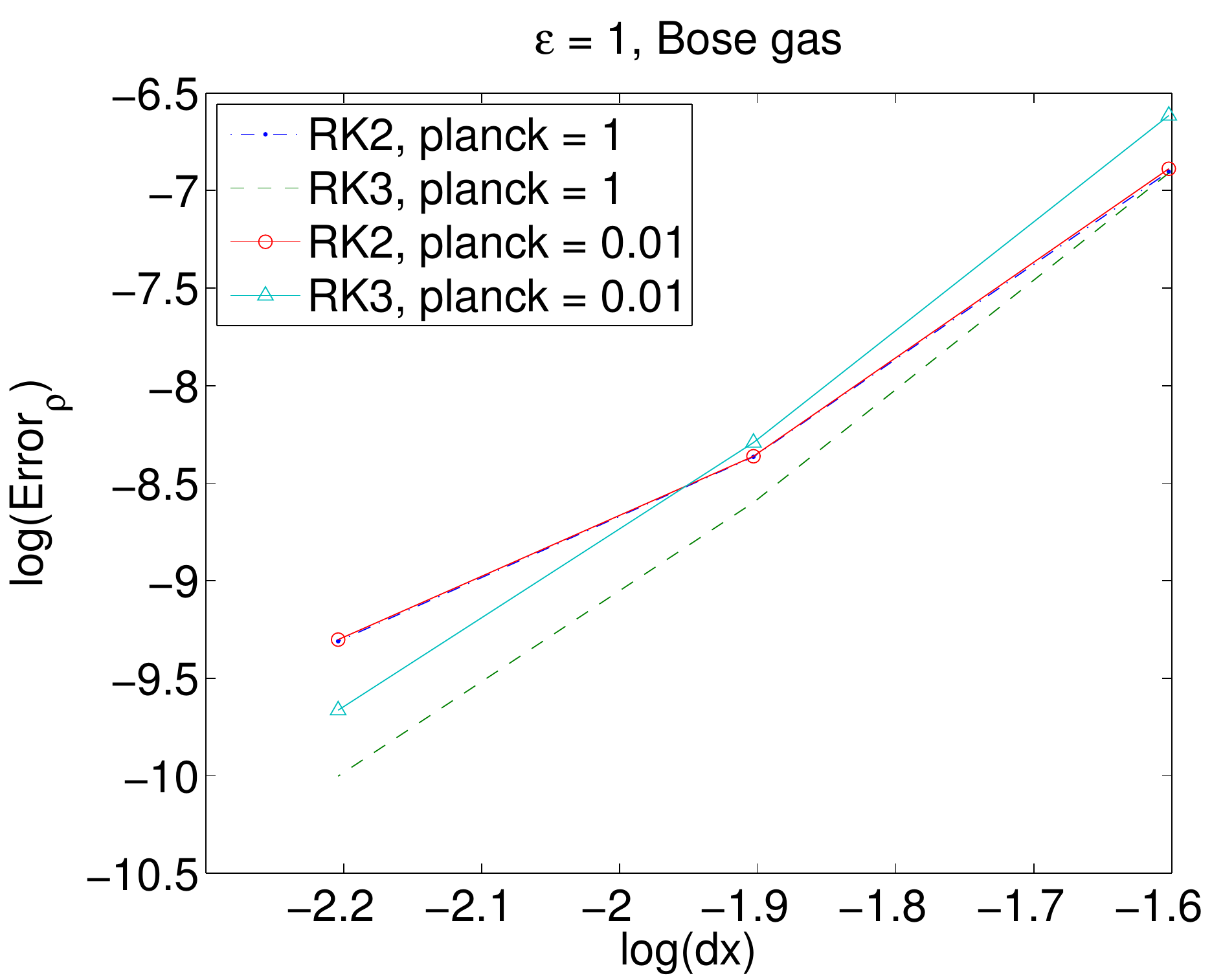}
   \includegraphics[width=3.2in]{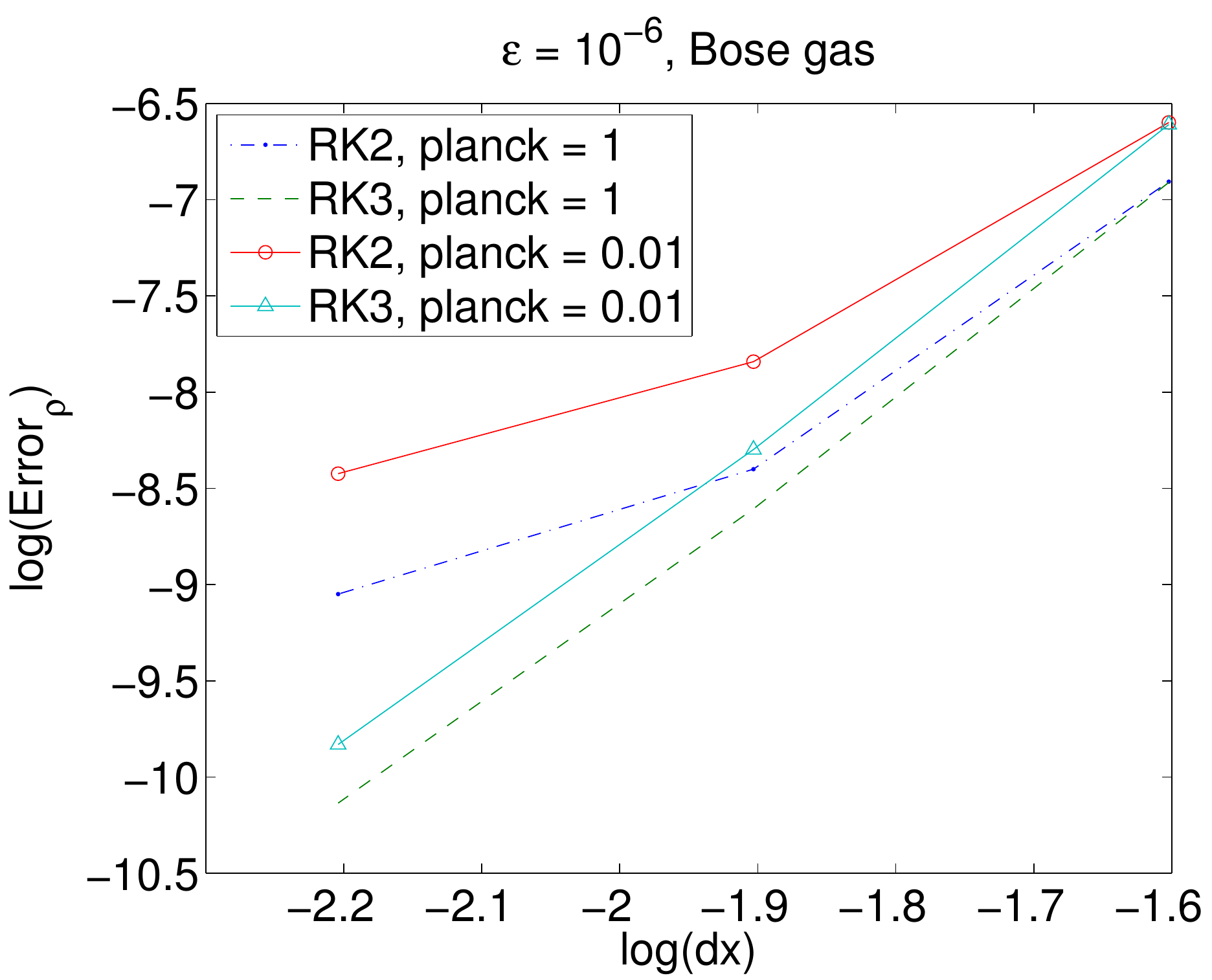}}
   \subfigure[Fermi gas]
   {\includegraphics[width=3.2in]{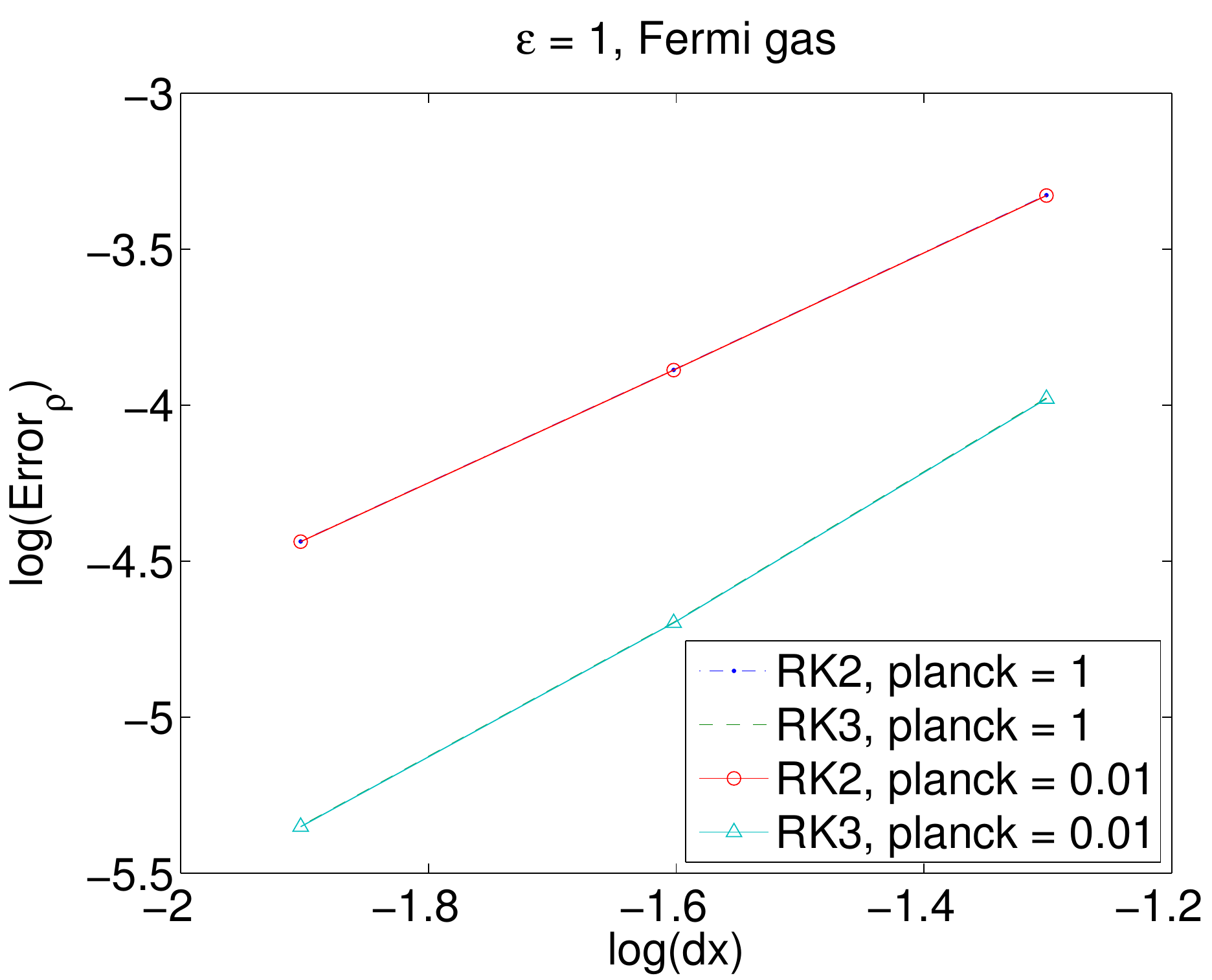}
   \includegraphics[width=3.2in]{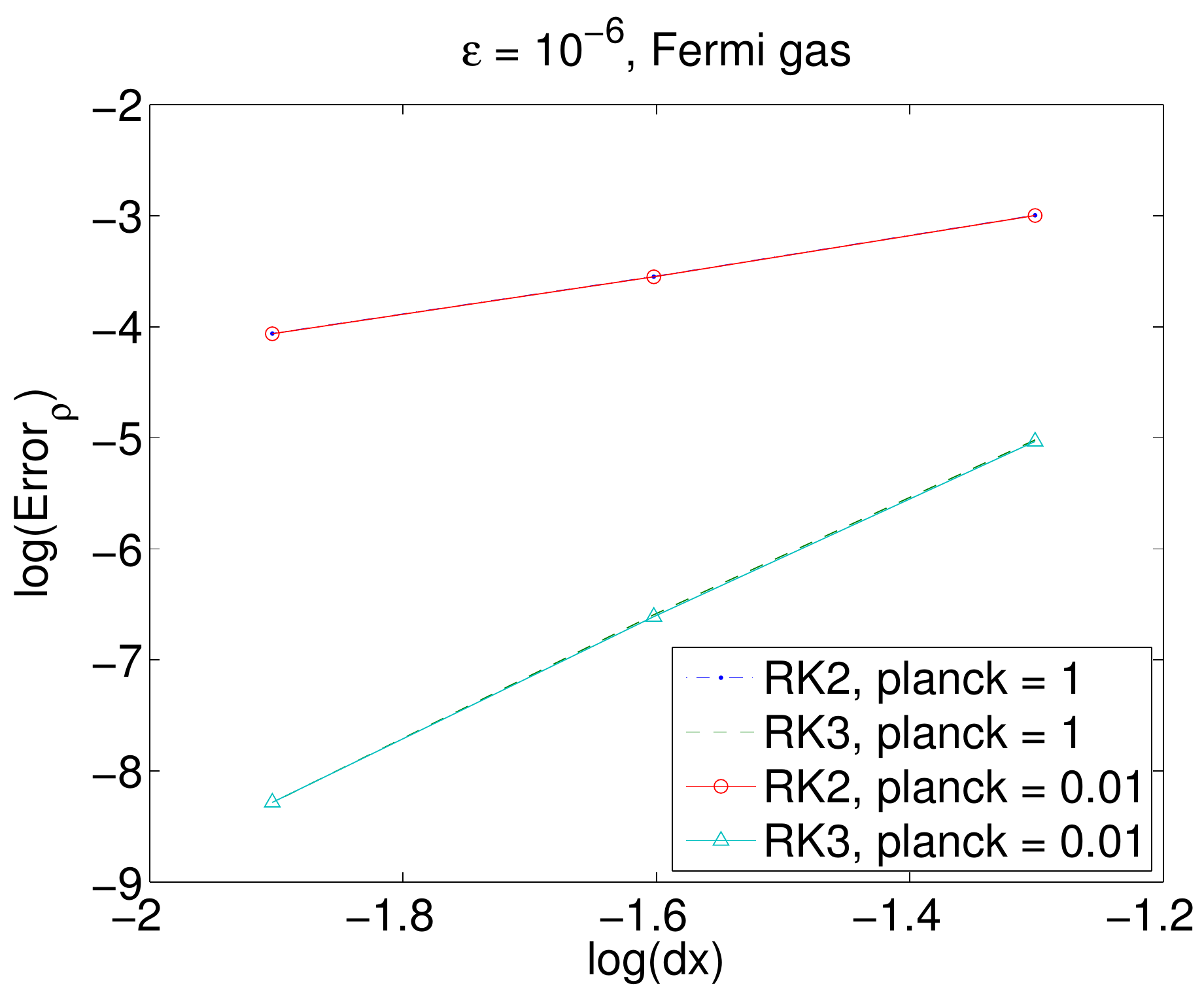}}
 \end{center}\caption{Convergence rate test. The left column is for 
 $\varepsilon=1$ and $\varepsilon=10^{-6}$ on the right. The top figures are 
 for Bose gas and the bottom ones are for Fermi gas.}\label{fig_conv_Bose}
\end{figure}

In each subfigure in Figure~\ref{fig_conv_Bose}, we show the convergence rate with $\theta_0=1$ 
and $\theta_0=10^{-2}$ using Exp-RK2 and Exp-RK3. We perform the same test for both 
Bose gas and Fermi gas in both kinetic regime and fluid regime. The numerical 
results are in good agreement with our theoretical expectation.

\section{Conclusions}
In this paper we have extended the numerical approach recently introduced in~\cite{LP_ExpRKinhomoBoltzmann} to the case of the quantum Boltzmann equation. In particular, we have shown how to derive high-order asymptotic preserving schemes which work uniformly with respect to the Planck constant. Numerical results for second and third order methods confirm the robustness and accuracy of the present method. We did not tackle the issue of the formation of the Bose-Einstein condensate since this involves also a careful choice of the velocity discretization whereas here we concentrate our attention on the time discretization problem only. In our future research we will focus on these challenging aspects.

\vspace{0.3in}
\noindent{\bf Acknowledgments.} 
We thank CSCAMM, University of Maryland for holding the conference ``Quantum 
Systems: A Mathematical Journey from Few to Many Particles'' in May 2013, at which this
work was initiated.

\appendix
\newcommand{\appsection}[1]{\let\oldthesection\thesection
\renewcommand{\thesection}{Appendix \oldthesection}
\section{#1}\let\thesection\oldthesection}

\appsection{Derivation of $\partial_t \mathcal{M}_q$}

In this appendix, we give the details of the derivation of (\ref{Mqt2}). Our goal is to represent $\partial_tz$ and $\partial_tT$ in equation (\ref{Mqt}) in terms of $\partial_t\rho$ and $\partial_te$. 

First, combining the two equations in system (\ref{22system}) gives
\begin{equation}
\frac{Q_{\frac{d}{2}}^{\frac{d}{2}+1}(z)}{Q_{\frac{d}{2}+1}^{\frac{d}{2}}(z)}=\theta_0 \left(\frac{d}{4\pi e}\right)^{\frac{d}{2}} \rho.
\end{equation}
Therefore, we define a function $F(z)$ such that
\begin{equation}
y=F(z)=\frac{Q_{\frac{d}{2}}^{\frac{d}{2}+1}(z)}{Q_{\frac{d}{2}+1}^{\frac{d}{2}}(z)},
\end{equation}
and a function $G(y)$ such that
\begin{equation}
z=G(y)=F^{-1}(y).
\end{equation}
Then we have
\begin{align} \label{G}
G'(y)=\frac{1}{F'(z)}=\frac{Q_{\frac{d}{2}+1}^d(z)}{\left(\frac{d}{2}+1\right)Q_{\frac{d}{2}}^{\frac{d}{2}}(z)Q_{\frac{d}{2}}'(z)Q_{\frac{d}{2}+1}^{\frac{d}{2}}(z)-\frac{d}{2}Q_{\frac{d}{2}+1}^{\frac{d}{2}-1}(z)Q_{\frac{d}{2}+1}'(z)Q_{\frac{d}{2}}^{\frac{d}{2}+1}(z)}.
\end{align}
For the Bose-Einstein/Fermi-Dirac function, one has the following nice property (see \cite{Pathria})
\begin{equation} \label{prop}
zQ_{\nu}'(z)=Q_{\nu-1}(z).
\end{equation}
Using (\ref{prop}) in (\ref{G}),
\begin{align}
G'(y)&=\frac{1}{F'(z)}=\frac{zQ_{\frac{d}{2}+1}^d(z)}{\left(\frac{d}{2}+1\right)Q_{\frac{d}{2}}^{\frac{d}{2}}(z)Q_{\frac{d}{2}-1}(z)Q_{\frac{d}{2}+1}^{\frac{d}{2}}(z)-\frac{d}{2}Q_{\frac{d}{2}+1}^{\frac{d}{2}-1}(z)Q_{\frac{d}{2}}(z)Q_{\frac{d}{2}}^{\frac{d}{2}+1}(z)} \nonumber\\
&=\frac{zQ_{\frac{d}{2}+1}^{\frac{d}{2}+1}(z)}{Q_{\frac{d}{2}}^{\frac{d}{2}}(z)\left[\left(\frac{d}{2}+1\right)Q_{\frac{d}{2}-1}(z)Q_{\frac{d}{2}+1}(z)-\frac{d}{2}Q_{\frac{d}{2}}^{2}(z)\right]}.
\end{align}
From the second equation of (\ref{22system}) we know
\begin{equation}
Q_{\frac{d}{2}+1}(z)=\frac{2e}{dT}Q_{\frac{d}{2}}(z),
\end{equation}
then
\begin{align}
G'(y)=\frac{z\left(\frac{2e}{dT}\right)^{\frac{d}{2}}}{\left(\frac{d}{2}+1\right)Q_{\frac{d}{2}-1}(z)-\frac{d^2T}{4e}Q_{\frac{d}{2}}(z)}.
\end{align}

Note that
\begin{align}
z=G\left(\theta_0 \left(\frac{d}{4\pi e}\right)^{\frac{d}{2}} \rho\right), \quad T=\frac{\theta_0^{\frac{2}{d}}}{2\pi}\left(\frac{\rho}{Q_{\frac{d}{2}}(z)}\right)^{\frac{2}{d}},
\end{align}
so we have
\begin{align}
\partial_t z&=G'\left(\theta_0 \left(\frac{d}{4\pi e}\right)^{\frac{d}{2}} \rho\right)\theta_0\left(\frac{d}{4\pi}\right)^{\frac{d}{2}}\left(\frac{1}{e^{\frac{d}{2}}}\partial_t\rho-\frac{d}{2}\frac{\rho}{e^{\frac{d}{2}+1}}\partial_te\right)\nonumber\\
&=\frac{zQ_{\frac{d}{2}}(z)}{\left(\frac{d}{2}+1\right)Q_{\frac{d}{2}-1}(z)-\frac{d^2T}{4e}Q_{\frac{d}{2}}(z)}\left(\frac{1}{\rho}\partial_t\rho-\frac{d}{2e}\partial_te\right),
\end{align}
and
\begin{align}
\partial_tT&=\frac{\theta_0^{\frac{2}{d}}}{\pi d}\left(\frac{\rho^{\frac{2}{d}-1}}{Q_{\frac{d}{2}}^\frac{2}{d}(z)}\partial_t\rho-\frac{\rho^{\frac{2}{d}}Q_{\frac{d}{2}}'(z)}{Q_{\frac{d}{2}}^{\frac{2}{d}+1}(z)}\partial_tz\right)=\frac{\theta_0^{\frac{2}{d}}}{\pi d}\left(\frac{\rho^{\frac{2}{d}-1}}{Q_{\frac{d}{2}}^\frac{2}{d}(z)}\partial_t\rho-\frac{\rho^{\frac{2}{d}}Q_{\frac{d}{2}-1}(z)}{zQ_{\frac{d}{2}}^{\frac{2}{d}+1}(z)}\partial_tz\right)\nonumber\\
&=\frac{2T}{d}\frac{1}{\rho}\rho_t-\frac{2T}{d}\frac{Q_{\frac{d}{2}-1}(z)}{Q_{\frac{d}{2}}(z)}\frac{1}{z}\partial_tz=\frac{2T}{d}\frac{1}{\rho}\rho_t-\frac{2T}{d}\frac{Q_{\frac{d}{2}-1}(z)}{\left(\frac{d}{2}+1\right)Q_{\frac{d}{2}-1}(z)-\frac{d^2T}{4e}Q_{\frac{d}{2}}(z)}\left(\frac{1}{\rho}\partial_t\rho-\frac{d}{2e}\partial_te\right).
\end{align}
Therefore,
\begin{align} 
&\frac{1}{z}\partial_tz+\frac{(v-u)^2}{2T^2}\partial_tT=\frac{Q_{\frac{d}{2}}(z)}{\left(\frac{d}{2}+1\right)Q_{\frac{d}{2}-1}(z)-\frac{d^2T}{4e}Q_{\frac{d}{2}}(z)}\left(\frac{1}{\rho}\partial_t\rho-\frac{d}{2e}e_t\right)  \nonumber\\
&+\frac{(v-u)^2}{dT}\frac{1}{\rho}\partial_t\rho-\frac{(v-u)^2}{dT}\frac{Q_{\frac{d}{2}-1}(z)}{\left(\frac{d}{2}+1\right)Q_{\frac{d}{2}-1}(z)-\frac{d^2T}{4e}Q_{\frac{d}{2}}(z)}\left(\frac{1}{\rho}\partial_t\rho-\frac{d}{2e}\partial_te\right)\nonumber\\
&=\left[\frac{Q_{\frac{d}{2}}(z)}{\left(\frac{d}{2}+1\right)Q_{\frac{d}{2}-1}(z)-\frac{d^2T}{4e}Q_{\frac{d}{2}}(z)}+ \frac{(v-u)^2}{dT}\left(1-\frac{Q_{\frac{d}{2}-1}(z)}{\left(\frac{d}{2}+1\right)Q_{\frac{d}{2}-1}(z)-\frac{d^2T}{4e}Q_{\frac{d}{2}}(z)}  \right)  \right]\frac{1}{\rho}\partial_t\rho \nonumber\\
&+\left[\frac{(v-u)^2}{2eT}\frac{Q_{\frac{d}{2}-1}(z)}{\left(\frac{d}{2}+1\right)Q_{\frac{d}{2}-1}(z)-\frac{d^2T}{4e}Q_{\frac{d}{2}}(z)} -\frac{d}{2e}\frac{Q_{\frac{d}{2}}(z)}{\left(\frac{d}{2}+1\right)Q_{\frac{d}{2}-1}(z)-\frac{d^2T}{4e}Q_{\frac{d}{2}}(z)}   \right]\partial_te.
 \end{align}
 Then if we define $M(z)$ and $N(z)$ as in (\ref{MN}), (\ref{Mqt2}) follows readily from the above equation.

\end{document}